\documentclass[12pt,centertags,oneside]{amsart}
\usepackage{amsmath,amstext,amsthm,amscd,typearea,hyperref,stmaryrd}
\usepackage{amssymb}
\usepackage{a4wide}
\usepackage[mathscr]{eucal}
\usepackage{mathrsfs}
\usepackage{typearea}
\usepackage{charter}
\usepackage{pdfsync}
\usepackage[a4paper,width=16.2cm,top=3cm,bottom=3cm]{geometry}
\usepackage{xcolor}

\numberwithin{equation}{section}

\newtheorem{theorem}{Theorem}[section]
\newtheorem{definition}[theorem]{Definition}
\newtheorem{proposition}[theorem]{Proposition}

\newtheorem{lemma}[theorem]{Lemma}

\newcommand{\dist}{{\rm dist}}

\renewcommand{\Im}{{\rm Im}}

\newcommand{\C}{\mathbb{C}}
\newcommand{\D}{\mathbb{D}}
\newcommand{\I}{\mathbb{I}}
\newcommand{\E}{\mathbb{E}}
\newcommand{\F}{\mathbb{F}}
\renewcommand{\H}{\mathbb{H}}
\newcommand{\M}{\mathbb{M}}

\newcommand{\Z}{\mathbb{Z}}

\newcommand{\R}{\mathbb{R}}
\newcommand{\T}{\mathbb{T}}

\newcommand{\G}{\mathbb{G}}
\newcommand{\U}{\mathbb{U}}

\renewcommand{\P}{\mathbb{P}}

\newcommand{\CommaBin}{\mathbin{\raisebox{0.5ex}{,}}}

\newcommand{\Bf}{\mathfrak{B}}
\newcommand{\Rf}{\mathfrak{R}}

\newcommand{\Zariski}{{\rm Zar}}

\title{Algebraic flows on commutative complex Lie groups}

\author{Tien-Cuong Dinh}
\address{Department of Mathematics, National University 
of Singapore, 10 Lower Kent Ridge Road, Singapore 119076. {\tt http://www.math.nus.edu.sg/$\sim$matdtc}}
\email{matdtc@nus.edu.sg}

\author{Duc-Viet Vu}

\address{University of Cologne, Mathematical Institute, Weyertal 86-90
D-50931 K\"oln, Germany.}
\email{vuduc@math.uni-koeln.de}

\begin{document}

\maketitle

\begin{abstract}
We recover results by  Ullmo-Yafaev and Peterzil-Starchenko on the closure of the image of an algebraic variety in a compact complex torus. Our approach uses directed closed currents and allows us to extend the result for dimension $1$ flows  to the setting of commutative complex Lie groups which are not necessarily compact.  A version of the classical Ax-Lindemann-Weierstrass theorem for commutative complex Lie groups is also given. 
\end{abstract}

\medskip\medskip

\noindent
{\bf MSC 2010:} 32H, 32U, 37xx

 \medskip

\noindent
{\bf Keywords:}  algebraic flow, directed closed current, Ax-Lindemann-Weierstrass theorem.

\section{Introduction}\label{s:intro}

Let $\E$ be a complex vector space of dimension $n$ and let $\Gamma$ be a discrete additive subgroup of $\E$. The quotient space $\T:=\E/\Gamma$ is a connected commutative complex Lie group. It is 
called {\it a (complex) semi-torus of dimension $n$} if $\Gamma$ spans $\E$ as a vector space over $\C$, see e.g. \cite{NW}.  Define  $\Gamma_\R:= \Gamma \otimes_\Z \R.$  When the rank of $\Gamma$ is maximal, i.e., equal to $2n$, we obtain a complex torus.  By a \emph{sub-semi-torus} (\emph{complex sub-torus}) of $\T,$ we mean a closed connected  complex Lie subgroup of $\T$ which is itself a semi-torus (a complex torus).  
Denote by $\pi:\E\to\T$ the canonical projection from $\E$ to $\T$. 

Let $X$ be an  irreducible complex algebraic subset  of $\E$. In \cite{UY2}, Ullmo-Yafaev  posed the problem of description of the usual topological closure of 
$\pi(X)$ when $\T$ is a torus. This problem is related to the classical Ax-Lindemann-Weierstrass (ALW for short) theorem that will be discussed later. Ullmo-Yafaev solved this problem in the case where $\dim X=1$. Our first main theorem is the following which extends their result to general commutative complex Lie groups, see also Theorems \ref{t:compact-1D} and \ref{t:1D-Gamma}.

\begin{theorem} \label{t:main-1}
Let $\T:=\E/\Gamma$ be a commutative complex Lie group as above and let $\pi:\E\to \T$ be the canonical projection. Let $X$ be an irreducible complex algebraic subset of dimension $1$ of $\E$. Then $\overline{\pi(X)}$ is the union of $\pi(X)$ and a finite number of translated connected closed real Lie subgroups of $\T.$
\end{theorem}

Note that one can compactify $X$ by adding finitely many points at infinity.  
Consider a germ of $X$ at a point  $x_\infty$ at infinity and $x$ in this germ going to $x_\infty$. Then
the set of cluster values of $\pi(x)$  is either empty or  a finite union of translated real Lie subgroups of $\T$.

When $\T$ is a torus, we get exactly one translated real sub-torus, denoted by $\T_j$, for each germ of $X$ at infinity. This is  the Ullmo-Yafaev's theorem mentioned above. 
Ullmo and Yafaev constructed measures 
which turn out to be the Haar measures on $\T_j$. Their study used techniques of oscillatory integrals. Here, we will construct sequences of directed closed currents of suitable dimensions whose supports are shown to be equal to $\T_j$. Our approach can be used for general commutative Lie groups which are not necessarily compact. Theorem \ref{t:main-1} above is a consequence of Theorems \ref{t:compact-1D} and \ref{t:1D-Gamma} in Sections \ref{s:compact-1D} and \ref{s:1D}.

We will show later at the end of Section \ref{s:higherdim} that the problem for algebraic flows of higher dimension for semi-tori  is of a different nature. In the case of tori, we obtain, by using similar techniques and other basic tools from complex geometry, a new proof to the following result due to  Peterzil-Starchenko \cite{PS}, see also \cite{PS2} and Theorem \ref{t:main_high_compact}.

\begin{theorem}[Peterzil-Starchenko] \label{t:main-2} 
Let  $\E, \T, \pi$ be as above with $\T$ a torus. Let $X$ be an irreducible complex algebraic subset  of $\E$.
Then there are finitely many  complex algebraic  subsets $C_1, \ldots,C_m$ of  $\E$, of dimension strictly smaller than $\dim X$, and  real sub-tori $\T_1, \ldots, \T_m$ of $\T$ such that 
$$\overline{\pi(X)} = \pi(X)\cup \bigcup_{j=1}^m \big(\pi(C_j)+ \T_j \big).$$
\end{theorem}

In their proof,  Peterzil and Starchenko used o-minimal theory which allows them to get a similar statement for categories of sets other than complex algebraic ones. 
Our approach will be presented in Section \ref{s:higherdim}. It uses the case of dimension 1 and is somehow more explicit.
Several steps of the proof  are quantifiable and this may be useful for  applications. 
Theorem \ref{t:main-2} will be obtained as a consequence of Theorem \ref{t:main_high_compact}.
When $\T$ is non-compact, the same statement is no more true in general as shows an example 
at the end of Section \ref{s:higherdim}. In this setting, 
we may expect a similar result with $\F_j$ real vector spaces and $C_i$ real (semi-)algebraic sets.

When $\T$ is a semi-abelian variety, the classical ALW theorem says that the Zariski closure of $\pi(X)$ is a translated sub-semi-abelian variety of $\T$, see Ax \cite{Ax} and Noguchi \cite{Noguchi}.  The classical ALW still holds for non-projective compact torus and in this case the Zariski  closure of $\pi(X)$ is a translated sub-torus, see  Pila-Zannier \cite{PZ} or Demailly \cite{Demailly_hyper}. We also refer to Paun-Sibony \cite{PaS} for a generalization of ALW theorem to the case of parabolic Riemann surfaces in a torus.
Using the technique of directed closed currents, we get the following version of  the ALW theorem for general commutative complex Lie groups, see also Theorem \ref{t:ALW-1D}.

\begin{theorem} \label{t:ALW} 
Let $\T=\E/\Gamma$ be a commutative complex Lie group as above and $X$ an  irreducible algebraic subset of $\E.$ 
Denote by $\overline{\pi(X)}^\Zariski$ the complex analytic Zariski closure of $\pi(X)$ in $\T$. There exist a sub-semi-torus $\T'$ of $\T$ and an irreducible analytic subset $S$ of $\T/\T'$ such that if   $p:\T\to \T/ \T'$ denotes the natural projection then
\begin{enumerate}
\item[(a)] the image of $X$ by  $p\circ\pi$ is a dense Zariski open subset of $S$;
\item[(b)]  we have  $\overline{\pi(X)}^\Zariski=p^{-1}(S).$
\end{enumerate}
Moreover, when $\pi(X)$ is relatively compact in $\T$,  then $p\big(\pi(X)\big)$ is a point and $\overline{\pi(X)}^\Zariski$ is a sub-semi-torus of  $\T$.
\end{theorem}

The proof of the last result will be given in  Section \ref{s:ALW}.

\medskip

\noindent
{\bf Notation.} 
Throughout the paper, $\E$ is a complex vector space of dimension $n$, $\T=\E/\Gamma$ is a commutative complex Lie group of dimension $n$ and $\pi:\E\to\T$ denotes the canonical projection. If $\F$ is a real or complex vector subspace of $\E$, denote by $\Pi_\F:\E\to\E/\F$ the canonical projection. If $f$ is a map with values in $\E$, then we define $f_\F:=\Pi_\F\circ f$. 

Fix a Hermitian metric on $\E$ associated to a K\"ahler $(1,1)$-form $\omega_\E$ with constant coefficients, or equivalently, invariant by translation. This metric induces a metric on $\T$. The form $\omega_\E$ also induces a K\"ahler $(1,1)$-form on $\T$ that we still denote by $\omega_\E$ for simplicity. We apply the same rule for all differential forms with constant coefficients on $\E$.

Define $\I:=(0,1)$, $\D^*:=\D\backslash\{0\}$ and $\overline \D^*:=\overline\D\backslash \{0\}$, where $\D:=\{|x|<1\}$ is the unit disc in $\C$. Define $\D(r):=\{|x|<r\}$ the disc of center 0 and radius $r$ in $\C$. 
 For any subset $\Theta$ of $\R/2\pi\Z$, denote by $\D_\Theta$ the union of 
the radii $L_\theta:=\{x=re^{i\theta},\ r\in(0,1)\}$ with $\theta\in\Theta$. Beside $x$, we also use some coordinate $x'=\lambda x+O(x^2)$ on a disc $\D'(0,\rho):=\{|x'|<\rho\}$,  and denote by $\D'_\Theta$ the union of the radii $L'_\theta:=\{x'=re^{i\theta},\ r\in(0,\rho)\}$ with $\theta\in\Theta$.

For two functions $g(x)$ and $h(x)$ in a neighborhood of $0$ in $\C$, we write $g(x)=\Theta(h(x))$ as $x$ tends to $0$ if we have both $g(x)=O(h(x))$ and $h(x)=O(g(x))$ when $x$ tends to $0$.

\medskip
\noindent
{\bf Acknowledgments.} 
The  first author is supported by Provost's Chair and Tier 1 Grants C-146-000-047-001 and  R-146-000-248-114  from National University of Singapore (NUS).  The second author would like to thank the Alexander von Humboldt Foundation for its financial support.
This paper was prepared in part during visits of the authors at Korea Institute for Advanced Study (KIAS), NUS and University of Cologne. 
They would like to thank these institutions for their hospitality and also Emmanuel Ullmo who introduced them to this research topic.

\section{One dimensional  flows with compact support} \label{s:compact-1D}

In this section, we will give some basic properties for currents on $\E$ or $\T$ and we will study flows of dimension 1 whose supports are contained in a compact subset of $\T$. Until Theorem \ref{t:compact-1D}, we don't need to assume this compactness property.

Let $V$ be a real affine subspace  of $\E$ of dimension $l\ge 1$.  A current $S$  on $\E$ (resp. $\T$) is said to be \emph{directed} by $V$ if we have 
$S \wedge \Phi=0$ for every $1$-form $\Phi$ with constant coefficients whose restriction to $V$ (resp. to the immersed manifold $\pi(V)$) vanishes.  
Let $\F$ be the vector subspace of $\E$ of dimension $l$ which is parallel to $V$. We also consider it as a group acting on both $\E$ and $\T$.

\begin{lemma} \label{l:directed-current}
If a closed current $S$ of dimension $l$ and of order $0$ on $\E$ (resp. $\T$) is directed by $V$, then it is invariant by $\F$.
\end{lemma}
\proof
It is enough to consider the case of a current on $\E$ since we can always pull-back currents from $\T$ to $\E$. 
Choose a real coordinate system $x=(x_1,\ldots,x_{2n})$ on $\E$ such that $\F$ is given by $x_{l+1}=\cdots=x_{2n}=0$. 
Since $S$ is directed by $\F$, we have $S\wedge dx_j=0$  for $j=l+1,\ldots,2n$. 
So we can write 
$S=h dx_{l+1} \wedge  \ldots \wedge dx_{2n}$,
where $h$ is a 0-current of order 0 on $\E$. Since $S$ is closed, we easily deduce that $\partial h/\partial x_j=0$  for $j=1,\ldots,l$. 

Therefore, if $\Pi_\F:\E\to \E/\F$ denotes the canonical projection, there is a 0-current $h'$ on $\E/\F$ such that $h=\Pi_\F^*(h')$.
Define $\nu:=h'dx_{l+1} \wedge  \ldots \wedge dx_{2n}.$
We have $S=\Pi_\F^*(\nu)$. Finally, $\nu$ is a measure since $S$ is of order 0. The lemma follows easily. 
\endproof

\begin{proposition} \label{p:Haar-current} 
There exists,  up to a multiplicative constant, a unique  non-zero closed current $R$ of dimension $l$ and of order $0$ on $\T$, directed by $V$, whose support is contained in $\overline{\pi(V)}$. Moreover, the support of such a current is equal to $\overline{\pi(V)}$.
\end{proposition} 

Note that $\overline{\pi(V)}$ is a  translated connected closed Lie subgroup of $\T$. So there is a real affine subspace $W$ of $\E$ such that $\pi(W)=\overline{\pi(V)}$. 

\begin{definition} \rm
We say that $R$ is {\it a Haar current} associated with the affine subspace $V$.
\end{definition}

We will see in the proof below that when $V$ is a complex affine space then Haar currents associated to $V$ are all positive or negative.

\proof[Proof of Proposition \ref{p:Haar-current}]
Let $\nu$ be a standard volume form on $\Pi_\F(W)$ and define $S:=\Pi_\F^*(\nu)$. Clearly, if an element of $\Gamma$ preserves $W$ then it preserves $S$. It follows that $S$ induces a closed current $R$ of dimension $l$ and of order 0 on $\T=\E/\Gamma$, directed by $V$, with support in $\pi(W)=\overline{\pi(V)}$. 

Let $R'$ be any closed current of dimension $l$  and of order 0 on $\T=\E/\Gamma$, directed by $V$, with support in $\pi(W)=\overline{\pi(V)}$. We need to show that $R'$ is proportional to $R$. Observe that $\pi^*(R')$ is a closed current of dimension $l$ and of order 0 on $\E$ with support in $W+\Gamma$ which is a union of affine spaces parallel to $W$. Denote by $S'$ its restriction to $W$ which is also a closed current of dimension $l$ and of order 0, directed by $V$, or equivalently, directed by $\F$. By Lemma \ref{l:directed-current} and its proof, 
there is a  locally finite measure $\nu'$ on $\E/\F$ such that $S'=\Pi_\F^*(\nu')$. 

Observe now that $\nu'$ is supported by $\Pi_\F(W)$. Let $\G$ denote the real vector subspace of $\E$ obtained from $W$ by a translation. Since $S'$ is invariant by translations by vectors in $\F$, the current $R'$ is invariant by the action of the group $\pi(\F)$ on $\T$. It follows that $R'$ is invariant by $\overline{\pi(\F)}$, or equivalently, $S'$ is invariant by translations by vectors in $\G$. We then deduce that the measure $\nu'$ on $\Pi_\F(W)$ is invariant by the action of $\G$ on $\E/\F$. Thus, it is proportional to $\nu$. The proposition follows.
\endproof

We now consider a situation slightly more general than the one in Theorem \ref{t:main-1}.
Let $f:\overline\D^*\to \E$ be a smooth map which is holomorphic in $\D^*$ and has a polynomial growth at $0$, that is, 
$\|f(x)\|=O(x^{-d})$ as $x\to 0$ for some integer $d$. Here, $x$ denotes the standard complex coordinate on $\C$. 
Our goal is to describe  the set $\Lambda_f$ of cluster values of $\pi\big( f(x)\big)$ in $\T$ when $x$ tends to $0$. 
Observe that when $d\leq 0$, the map $f$ is bounded and can be extended holomorphically through 0. It follows that $\Lambda_f=\pi(f(0))$ in this case.
So, from now on, we assume the following property.

\medskip\noindent
{\bf (H0)} The map $f$ is not bounded.

\medskip\noindent
So there is a minimal integer $d\geq 1$ such that
$$\|f(x)\|= \Theta(x^{-d}).$$

Recall that $\Pi_\F:\E\to\E/\F$ denotes the canonical projection 
for any vector subspace $\F$ of $\E$, and we define $f_\F:=\Pi_\F\circ f$.

\begin{lemma} \label{l:F-l}
There are a unique integer $1\leq k\leq n$ and  a unique sequence of integers $d_1>\cdots>d_{k+1}$ with $d_1=d$, $d_{k+1}=0$,  such that  
$$f(x)=\big(1+ O(x)\big)x^{-d_1}v_1+\cdots + \big(1+ O(x)\big) x^{-d_k} v_k + v_{k+1} + O(x) \quad \text{as} \quad x\to 0$$
for some vectors $v_1,\ldots,v_{k+1}$ in $\E$ with $v_1,\ldots,v_k$ $\C$-linearly independent.
Denote by $\F_l$ the complex vector subspace of dimension $l$ of $\E$, spanned by $v_1,\ldots,v_l$. Then 
$v_1$ is unique and $v_{l+1}$ is unique modulo $\F_l$. Moreover, we also have 
$\|f_{\F_l}\|=\Theta(x^{-{d_{l+1}}})$ for $1\leq l\leq k-1$ and $\|f_{\F_k}\|=O(1)$, 
and $\Lambda_f\subset \overline{\pi(V)}$ 
with  $V:=\F_k + v_{k+1}$. 
\end{lemma}
\proof
Clearly, we need to take $d_1:=d$ and 
$$v_1:=\lim_{x\to 0} x^{d_1} f(x).$$
So $v_1$ and $d_1$ are unique. By induction,  for $l\geq 0$ with $d_l\geq 1$, either we have $d_{l+1}=0$ and $\|f_{\F_l}\|=O(1)$ or 
$d_{l+1}$ should be the positive integer such that $\|f_{\F_l}\|=\Theta(x^{-d_{l+1}})$.
Then, we define 
$$v^0_{l+1}:=\lim_{x\to 0} x^{d_{l+1}} f_{\F_l}(x).$$

It is not difficult to see that this vector exists in $\E/\F_l$, and it is unique and non-zero.
In order to have the identity for $f(x)$ in the lemma, we should choose a vector $v_{l+1}\in \E$ such that $\Pi_{\F_l}(v_{l+1})=v^0_{l+1}$. 
Therefore, the choice of $v_{l+1}$ is unique modulo $\F_l$. We then 
define $\F_{l+1}$ as the complex vector space spanned by $\F_l$ and $v_{l+1}$. We end the inductive construction at the step $k$ when we get $d_{k+1}=0$. 
It is then easy to check that the obtained $d_l, v_l$ and $\F_l$ satisfy the  lemma, except for the last assertion.

Finally, when $x$ tends to 0, we have that $f_{\F_k}(x)\to  v^0_{k+1}$. Thus, the distance between $f(x)$ and $V$ tends to 0. 
The last assertion in the lemma follows. Note that the spaces $\F_l$ are canonically associated to the map $f$. In particular, $\F_k$ is the smallest vector subspace of $\E$ such that $f_{\F_k}$ can be extended to a holomorphic map from $\D$ to $\E/\F_k$.
\endproof

Consider now a  domain $\U$ with piecewise smooth boundary which is relatively compact in $\D^*$. Define $\U_a:= a \U$ for $a\in\C$ with $|a|\leq 1$. We have the following lemma. 

\begin{lemma} \label{l:mass}
Let  $\|v_1\|$ be the norm of $v_1$ with respect to $\omega_\E$ and define 
$$\lambda_0:=d^2 \|v_1\|^2\int_{x\in \U} |x|^{-2d-2} (idx\wedge d\overline x).$$
Then, the mass of the positive measure $f_* [\U_a]\wedge \omega_\E$ is equal to $\lambda_0 |a|^{-2d} + O(a^{-2d+1})$ when $a$ tends to $0$.
\end{lemma}
\proof  
Recall that $x$ is the standard complex coordinate on $\C$ and denote by $v:=\partial/\partial x$ the unit complex tangent vector at a point $x$. 
Denote by $df(x)$ the differential of $f$ at the point $x$ and define $u:=df(x)(v)$. The value of $\omega_\E$ at the point $f(x)$ is a quadratic form on the complex tangent space of $\E$ at this point and $\omega_\E(u,\overline u)$ denotes its values at $(u,\overline u)$. 
Then, using the expansion of $f$ in Lemma \ref{l:F-l}, the mass of $f_* [\U_a]\wedge \omega_\E$ is equal to
\begin{eqnarray*} 
\|f_* [\U_a]\wedge \omega_\E \| & = &  \int_{\U_a} f^* (\omega_\E)= \int_{x\in \U_a} \omega_\E\big(u, \overline u\big) dx\wedge d\overline x \\
& = &  d^2 \|v_1\|^2  \int_{x\in \U_a} \big(|x|^{-2d-2}+O(x^{-2d-1})\big) (i dx \wedge d \overline x) \\
& = & \Big[d^2 \|v_1\|^2\int_{x\in \U} |x|^{-2d-2} (idx\wedge d\overline x)\Big] |a|^{-2d} + O(a^{-2d+1}),
\end{eqnarray*}
where we used the change of variable $x\mapsto ax$. The lemma follows.
\endproof

Define 
$$\mu_{a}:=  \lambda_0^{-1} |a|^{2d}\pi_*\big(  f_* [\U_a]  \wedge \omega_\E \big),$$
where $\lambda_0$ is the constant in Lemma \ref{l:mass}. This is a positive measure of mass $1+O(a)$ on 
 $\T.$  

\begin{theorem} \label{t:compact-1D}  
We use the notation introduced above and assume that the hypothesis (H0) holds; in particular, the affine space $V$ is defined in Lemma \ref{l:F-l}.
Assume moreover that $\F_k$ is contained in the real vector space $\Gamma_\R$ spanned by $\Gamma$, or equivalently, $\overline{\pi(V)}$ is a translated real compact sub-torus of $\T$.
Let $\mu_V$ be the probability Haar  measure on  $\overline{\pi(V)}$.   Then $\mu_{a}$ converges to $\mu_V$ as $a$ goes to $0.$ Moreover, 
we have $\Lambda_f= \overline{\pi(V)}$ and  $\overline{\pi(f(\D^*))}=\pi(f(\overline\D^*))\cup  \overline{\pi(V)}$.
\end{theorem}

Until the end of this section, we work under the hypothesis of the above theorem. 
Since $\overline{\pi(V)}$ is compact, by  Lemma \ref{l:F-l},  the image of $\pi\circ f$ is contained in a compact subset of $\T$.
The first assertion of  Theorem \ref{t:compact-1D} implies that $\overline{\pi(V)} \subset \Lambda_f$. This, together with  Lemma \ref{l:F-l},  implies that $\Lambda_f= \overline{\pi(V)}$. Therefore, in order to prove Theorem \ref{t:compact-1D}, we only need to show that   $\mu_a$ converges to $\mu_V$. For this purpose, we will construct positive closed currents of suitable dimensions whose supports are equal to $\overline{\pi(V)}$.
They are associated to the  stratification $\F_1\subset \cdots \subset \F_k\subset \E$.

 For $1\leq l\leq k$, define the map $\tau_{l,a}: \overline \U\times \overline{\D}^{l-1} \rightarrow \E$  by 
$$\tau_{l,a}(x,t_1, \ldots, t_{l-1}) := f(ax)+ \sum_{j=1}^{l-1} a^{-d_j} t_jv_j .$$
Consider the positive current 
$S_{l,a}$ of bi-dimension $(l,l)$ on $\E$ defined by 
$$S_{l,a}:=c_{l,a}^{-1}(\tau_{l,a})_*\big(\Psi_l [\U\times\D^{l-1}]\big),$$
where $[\U\times\D^{l-1}]$ is the current of integration on $\U\times\D^{l-1}$ and 
$$\Psi_l(x,t_1,\ldots,t_{l-1}):=|x|^{-2(d-d_l)} \qquad \text{and} \qquad c_{l,a}:=  (2\pi)^{l-1}d_l^2|a|^{-2d_1-\cdots-2d_l}.$$
Equivalently, for every smooth $(l,l)$-form $\Phi$ with compact support on $\E$, we have 
$$\big\langle S_{l, a}, \Phi \big\rangle := c^{-1}_{l,a}  \int_{\U \times \D^{l-1}}  |x|^{-2(d-d_l)} (\tau_{l,a})^*(\Phi).$$ 
Denote for simplicity $t:=(t_1,\ldots,t_{l-1})$.

\begin{lemma} \label{l:vector}
Let  $u,u_1,\ldots,u_{l-1}$ be the images of the tangent vectors
$${\partial\over\partial x} \CommaBin \  {\partial\over\partial t_1} \ \CommaBin \  \cdots \  \CommaBin \  {\partial\over\partial t_{l-1}} $$
by the differential $d\tau_{l,a}(x,t)$ of $\tau_{l,a}$ at the point $(x,t)$. Then we have 
$$u\wedge u_1\wedge \ldots\wedge u_{l-1} = (-1)^ld_l a^{-d_1-\cdots-d_l}  x^{-d_l-1} v_1\wedge\ldots\wedge v_l +O(a^{-d_1-\cdots-d_l+1}).$$
Moreover, if $w_1,\ldots,w_{2l-1}$ belong to the family of vectors $u,\overline u,u_1,\overline u_1,\ldots,u_{l-1}, \overline u_{l-1}$, then
$$w_1\wedge \ldots\wedge w_{2l-1} = O(a^{-2d_1-\cdots-2d_l+1}).$$
\end{lemma}
\proof
Using the definition of $\tau_{l,a}$, we have $u_j=a^{-d_j} v_j$. Therefore, since $v_j\wedge v_j=0$,  the components of $u$ involving 
$v_1,\ldots,v_{l-1}$ do not contribute to the wedge product $u\wedge u_1\wedge\ldots\wedge u_{l-1}$. 
On the other hand, using the expansion of $f(x)$ in Lemma \ref{l:F-l}, we have 
$u=-d_la^{-d_l}x^{-d_l-1}v_l+O(a^{-d_l+1})$ modulo $v_1,\ldots, v_{l-1}$.
This implies the first assertion in the lemma.

Consider the second assertion. If the $w_j$ are not distinct, then their wedge-product vanishes. So we only need to consider a family $w_1,\ldots, w_{2l-1}$ which is obtained from 
the family $u,\overline u,u_1,\overline u_1,\ldots,u_{l-1}, \overline u_{l-1}$ by removing an element. If the removed element is $u$ or $\overline u$, the same argument as above gives the result. Otherwise, we assume that the removed element is $u_{j_0}$ for some $j_0$; the case of $\overline u_{j_0}$ can be obtained in the same way. As above, we have 
$$\overline u\wedge \overline u_1\wedge \ldots\wedge \overline u_{l-1} = O(a^{-d_1-\cdots-d_l})$$
and we need to estimate
$$u\wedge u_1\wedge\ldots\wedge u_{j_0-1}\wedge u_{j_0+1}\wedge \ldots\wedge u_{l-1}.$$
Using again that $u_j=a^{-d_j} v_j$ and $$u=-d_{j_0}a^{-d_{j_0}}x^{-d_{j_0}-1}v_{j_0}+O(a^{-d_{j_0}+1})$$ modulo $v_1,\ldots,v_{j_0-1},$   we see that the last wedge-product is of norm $O(a^{-d_1-\cdots-d_{l-1}})$.
The second assertion of the lemma follows easily.
\endproof

\begin{lemma} \label{l:mass-S-l}  
When $a$ tends to $0$, we have 
$$\|S_{l,a}\| = \|v_1\wedge\ldots\wedge v_l\|^2 \int_\U  |x|^{-2d-2}idx\wedge d\overline x + O(a) \qquad \text{and} \qquad \|dS_{l,a}\| =  O(a),$$
where the norm of the vector  $v_1\wedge\ldots\wedge v_l$ is with respect to the metric induced by  $\omega_\E$.
\end{lemma} 
\proof
Write for simplicity  $idt\wedge d\overline t:= (idt_1\wedge d\overline t_1)\wedge \ldots\wedge (idt_{l-1}\wedge d\overline t_{l-1})$. 
If $A$ is a point in $\U\times \D^{l-1}$ and $W$ is a tangent $q$-vector at $A$, denote by $\delta_A\otimes W$ the current of dimension $q$ such that 
$$\langle \delta_A\otimes W, \alpha\rangle := \langle W, \alpha(A)\rangle$$
for every smooth test form $\alpha$ of degree $q$. The last pairing is the number obtained using the contraction operator for the vectors $W$ and $\alpha(A)$.
We will consider disintegration of currents into currents of type $\delta_A\otimes W$.

Write
$$W:=\Big(i{\partial\over\partial x}\wedge {\partial\over\partial \overline x}\Big) \wedge \Big(i{\partial\over\partial t_1}\wedge {\partial\over\partial \overline t_1}\Big) \wedge \ldots \wedge \Big(i{\partial\over\partial t_{l-1}}\wedge {\partial\over\partial \overline t_{l-1}}\Big)$$
and 
$$W':=(iu\wedge \overline u)\wedge (iu_1\wedge \overline u_1)\wedge\ldots \wedge (iu_{l-1}\wedge \overline u_{l-1}).$$
By Lemma \ref{l:vector}, we have 
$$W' =  (2\pi)^{-l+1} c_{l,a}  |x|^{-2d_l-2} W'_l  + c_{l,a} O(a),$$
where
$$W'_l:= (iv_1\wedge \overline v_1)\wedge\ldots \wedge (iv_l\wedge \overline v_l).$$

Let $\mu_l$ denote the positive measure on $\U\times\D^{l-1}$ associated with the top degree form $|x|^{-2(d-d_l)} (idx\wedge d\overline x)\wedge (idt\wedge d\overline t)$.  Define $A':=\tau_{l,a}(A)$. By definition of $S_{l,a}$,  we have 
\begin{eqnarray*}
S_{l,a} & = &  c_{l,a}^{-1} (\tau_{l,a})_*\Big(\int_{A\in \U\times\D^{l-1}} (\delta_A\otimes W) d\mu_l(A)\Big) \ = \ c_{l,a}^{-1}
\int_{A\in \U\times \D^{l-1}} (\delta_{A'}\otimes W') d\mu_l(A) \\
& = &  (2\pi)^{-l+1} \int_{A\in \U\times \D^{l-1}} (\delta_{A'}\otimes W'_l) |x|^{-2d_l-2} d\mu_l(A) +O(a) 
\  = \ \widetilde\mu_l\otimes W'_l +O(a),
\end{eqnarray*}
where $\widetilde\mu_l$ is the positive measure defined by 
$$\widetilde\mu_l:= (2\pi)^{-l+1}  \int_{A\in \U\times \D^{l-1}} \delta_{A'} |x|^{-2d_l-2} d\mu_l(A).$$

Therefore, using the definition of $\mu_l$, we obtain that the mass $\|S_{l,a}\|$ of $S_{l, a}$ is equal to 
\begin{eqnarray*}
\|\widetilde\mu_l\| \|W'_l\| +O(a) &=& (2\pi)^{-l+1} \|v_1\wedge \ldots \wedge v_l\|^2  \int_{\U\times \D^{l-1}} |x|^{-2d_l-2} d\mu_l(A) +O(a) \\
& = & \|v_1\wedge \ldots \wedge v_l\|^2  \int_\U |x|^{-2d-2} idx\wedge d\overline x +O(a).
\end{eqnarray*}
This gives the first identity in the lemma.

We prove now the second identity. We have 
$$dS_{l,a}= c_{l,a}^{-1}(\tau_{l,a})_* \Big(d\big(\Psi_l [\U\wedge \D^{l-1}]\big)\Big)$$
and
$$d\big(\Psi_l [\U\wedge \D^{l-1}]\big) = \Psi_l [b(\U\times\D^{l-1})]+d\Psi_l \wedge [\U\times\D^{l-1}],$$
where $b(\U\times\D^{l-1})$ is the boundary of $\U\times \D^{l-1}$ with a suitable orientation. Then, the current $d\big(\Psi_l [\U\wedge \D^{l-1}]\big)$ can be obtained as an average of currents of the form 
$\delta_A\otimes \widetilde W$,
where $\delta_A$ is the Dirac mass at a point $A$ in $\overline\U\times\overline\D^{l-1}$ and $\widetilde W$ is a $(2l-1)$-tangent vector at $A$ with a bounded norm.
We deduce that $dS_{l,a}$ is an average of currents of the form $c_{l,a}^{-1}\delta_{A'}\otimes \widetilde W'$, where $\widetilde W'$ is the image of $\widetilde W$ by $d\tau_{l,a}(A)$. By Lemma \ref{l:vector}, the norm of $\widetilde W'$ is bounded by a constant times $|a|^{-2d_1-\cdots-2d_l+1}$. 
This, together with the definition of $c_{l,a}$, imply  the desired estimate.
\endproof

\begin{lemma} \label{l:limit-value-S-l}  
Any limit value of $\pi_*(S_{l,a})$, when $a$ tends to $0$, is a positive closed current of bi-dimension $(l,l)$ and of mass 
$$\lambda_l:= \|v_1\wedge\ldots\wedge v_l\|^2 \int_\U |x|^{-2d-2}idx\wedge d\overline x $$
which is invariant by $\F_l$.
\end{lemma}
\proof 
Observe that since the distances between the images of $\tau_{l,a}$ and $\F_k$ are bounded, the support of $\pi_*(S_{l,a})$ is contained in a fixed compact subset of $\T$. 
Consider a limit value $R_l$ of $\pi_*(S_{l,a})$.
By Lemma \ref{l:mass-S-l}, $R_l$ is a positive closed current of mass $\lambda_l$. It remains to check that $R_l$ is invariant by $\F_l$. 

By Lemma \ref{l:directed-current}, we only need to prove that $R_l$ is directed by $\F_l$. 
Let $\Phi$ be any 1-form with constant coefficients vanishing at $v_1,\overline v_1,\ldots, v_l, \overline v_l$. 
We have to show that  $S_{l,a}\wedge \Phi$  tends to 0 as $a$ tends to 0. 
In the proof of Lemma \ref{l:mass-S-l}, we obtained that $S_{l,a}=\widetilde\mu_l\otimes W'_l+O(a)$. Then, by definition of $W'_l$, we get 
$S_{l,a}\wedge\Phi=O(a)$. The lemma follows.
\endproof

\begin{proposition} \label{p:limit-S-k}  
When $a$ tends to $0$, the current $\pi_*(S_{k,a})$ tends to the Haar current of mass $\lambda_k$ associated with $V$ that we will denote by 
$R_k$.
\end{proposition}
\proof
Consider a limit value $R$ of $\pi_*(S_{k,a})$ when $a$ tends to 0. 
By Lemma \ref{l:limit-value-S-l}, this is a positive closed current of bi-dimension $(k,k)$ and of mass $\lambda_k$ directed by $\F_k$, or equivalently, by $V$. According to Proposition \ref{p:Haar-current}, it is enough to show that $R$ has support in the translated torus $\overline{\pi(V)}$.

Observe that the support of $S_{k,a}$ is contained in the image of $\overline \U\times\overline \D^{k-1}$ by $\tau_k$ that we denote by $\Sigma_{k,a}$. 
By definition of $\tau_k$, Lemma \ref{l:F-l}, and using that $\F_k\subset \Gamma_\R$, we see that the image of $\Sigma_{k,a}$ by $\Pi_{\F_k}$ tends to the image of $v_{k+1}$ in $\E/ \F_k$, or equivalently, we have 
$$\max_{z\in\Sigma_{l,a}}\dist(z,V) \to 0$$
as $a$ tends to 0. It follows that 
$$\max_{z\in \pi(\Sigma_{k,a})}\dist(z,\pi(V)) \to 0$$
as $a$ tends to 0. Since $\pi_*(S_{k,a})$ has support in $\pi(\Sigma_{k,a})$, we conclude that $R$ has support in the translated torus $\overline{\pi(V)}$. This ends the proof of the proposition.
\endproof

We choose differential $(1,0)$-forms $\Phi_l$ with constant coefficients on $\E$ for $l=1,\ldots,k$ such that we have  $\Phi_l(v_j)=1$ if $l=j$ and  $\Phi_l(v_j)=0$ if $l\not=j$ and $\Phi_l(\overline v_j)=0$ for every $1 \le j \le k$. If we consider a coordinate system $z=(z_1,\ldots,z_n)$ on $\E$ such that $v_j=\partial/\partial z_j$, then we can take $\Phi_j=dz_j$. 

\begin{proposition} \label{p:limit-S-l}
For every $1\leq l\leq k-1$, the current
$\pi_*(S_{l, a})$ tends to the positive closed current 
$$R_l:= R_k  \wedge (i \Phi_{l+1} \wedge \overline \Phi_{l+1})\wedge \ldots \wedge (i \Phi_k \wedge \overline \Phi_k).$$ 
Moreover, the support of $R_l$ is equal to the translated torus $\overline{\pi(V)}.$  
\end{proposition}
\proof
By Proposition \ref{p:limit-S-k}, $R_k$ is a Haar current associated with $V$. It is described in Lemma \ref{l:directed-current} and Proposition \ref{p:Haar-current}.
We easily see that the second assertion in the proposition is a consequence of the first one. We prove now the first assertion by induction on $l$.

Assume this property for $l+1,l+2,\ldots,k-1$. We only need to check that $\pi_*(S_{l, a})$ tends to $R_{l+1} \wedge (i \Phi_{l+1} \wedge \overline \Phi_{l+1})$. The hypothesis of induction implies that $\pi_*(S_{l+1, a})$ tends to $R_{l+1}$. 
For simplicity, we will only consider $a$ in an arbitrary  sequence converging to 0 such that $\pi_*(S_{l, a})$ tends to some current $R$. Then, we have to check that 
$R=R_{l+1} \wedge (i \Phi_{l+1} \wedge \overline \Phi_{l+1})$. By Lemma \ref{l:limit-value-S-l}, the current $R$ is invariant by $\F_l$.

Define 
$$W_l':=(iv_1\wedge \overline v_1)\wedge\ldots \wedge (iv_l\wedge \overline v_l) \qquad \text{and} \qquad W_{l+1}':= W_l' \wedge (iv_{l+1}\wedge \overline v_{l+1}).$$
In the proof of Lemma \ref{l:mass-S-l}, we proved that $S_{l,a}$ is equal, modulo a current of mass $O(a)$, to
$$\widetilde S_{l,a}:=\int_{A\in\U\times\D^{l-1}} (\delta_{A'}\otimes W_l') (2\pi)^{-l+1}|x|^{-2d_l-2}d\mu_l(A).$$
So the current $\pi_*(\widetilde S_{l,a})$ converges to $R$. 

Write $B:=(A,t_l)\in \U\times\D^l$ and $B':=\tau_{l+1}(B)=A'+a^{-d_l}t_lv_l$. We also have that $S_{l+1,a}$ is equal, modulo a current of mass $O(a)$, to
$$\widetilde S_{l+1,a}:=\int_{B\in\U\times\D^l} (\delta_{B'}\otimes W_{l+1}') (2\pi)^{-l}|x|^{-2d_{l+1}-2}d\mu_{l+1}(B).$$
The current $\pi_*(\widetilde S_{l+1,a})$ converges to $R_{l+1}$ by the hypothesis of induction. 

Observe  that the current
$\widetilde S_{l+1,a}\wedge (i\Phi_{l+1}\wedge \overline \Phi_{l+1})$ is equal to
\begin{eqnarray*}
\lefteqn{ \int_{B\in\U\times\D^l} (\delta_{B'}\otimes W_l') (2\pi)^{-l}|x|^{-2d_{l+1}-2} d\mu_{l+1}(B) } \\
& = & (2\pi)^{-1} \int_{t_l\in\D} \Big[\int_{A\in\U\times \D^{l-1}} (\delta_{B'}\otimes W_l') (2\pi)^{-l+1}|x|^{-2d_l-2}d\mu_l(A)\Big] (idt_l\wedge d\overline t_l) \\
& = & (2\pi)^{-1} \int_{t_l\in\D} \widetilde S_{l,a,t_l} (idt_l\wedge d\overline t_l),
\end{eqnarray*}
where $\widetilde S_{l,a,t_l}$ is the direct image of $\widetilde S_{l,a}$ by the translation by vector $a^{-d_l}t_lv_l$ which sends $A'$ to $B'$. 
Then, by taking the direct image by $\pi$ and the  limit when $a$ tends to 0, we have 
$$R_{l+1} \wedge (i \Phi_{l+1}\wedge \overline \Phi_{l+1}) = \lim_{a\to 0} (2\pi)^{-1} \int_{t_l\in\D} \pi_*(\widetilde S_{l,a,t_l}) (idt_l\wedge d\overline t_l).$$

Finally, recall that the limit $R$ of $\pi_*(\widetilde S_{l,a})$ is invariant by $\F_l$, in particular, it is invariant by the translation by the vector $a^{-d} t_l v_l$. It follows that the current $\pi_*(\widetilde S_{l,a,t_l})$ also converges to $R$. Hence, the last limit is equal to
$$(2\pi)^{-1}\int_{t_l\in\D} R (idt_l\wedge d\overline t_l)=R.$$
This completes the proof of the proposition.
\endproof

\proof[End of the proof of Theorem \ref{t:compact-1D}] 
Note that $d_1=d$ and $\tau_1(x)=f(ax)$. We have 
$$\mu_a=\lambda_0^{-1}|a|^{2d}\pi_*  f_* [\U_a] \wedge\omega_\E=\lambda_0^{-1}\pi_*(S_{1,a})\wedge\omega_\E,$$ 
where $\lambda_0>0$ is the constant given in Lemma \ref{l:mass}. It follows from Proposition \ref{p:limit-S-l} that $\mu_a$ converges to a constant times the 
measure $R_1\wedge \omega_\E$. The last measure is equal to the wedge-product of the Haar current associated with $V$ with a differential form with constant coefficients. Therefore, it should be proportional to the probability Haar measure on  $\overline{\pi(V)}$. 

Since $\mu_a$ is a positive measure of mass $1+O(a)$ with support in a fixed compact subset of $\T$, any limit value of $\mu_a$ is a probability measure. We conclude that $\mu_a$ converges to the probability Haar measure  on  $\overline{\pi(V)}$. As mentioned just after Theorem \ref{t:compact-1D}, the last property also implies the second assertion in that theorem and the proof is now complete.
\endproof

\section{One dimensional flows with non-compact support} \label{s:1D}

The goal of this section  is to  establish a convergence result similar to Theorem \ref{t:compact-1D} for $1$-dimensional  flows with non-compact support.  
We will use the notation introduced in Section \ref{s:compact-1D}. In this section, we assume the following property, see Theorem \ref{t:compact-1D}.

\medskip\noindent
{\bf (H1)}  $\F_k$ is not contained in the real vector space $\Gamma_\R$ spanned by 
$\Gamma$, or equivalently, the set $\overline{\pi(V)}$ is not compact in $\T$. 

\medskip

Let $1 \le \kappa \le k$  be the smallest integer such that $\C v_{\kappa} \not \subset \Gamma_\R$. 
Consider a radius $L_\theta:=\{x=re^{i\theta} \text{ with } r\in (0,1)\}$ for some $\theta\in\R/2\pi\Z$. 
We say that $L_{\theta}$ is an {\it almost $\Gamma$-radius} if $x^{-d_\kappa}v_\kappa$ belongs to $\Gamma_\R$ for $x\in L_{\theta}$, or equivalently, $e^{-id_\kappa\theta }v_\kappa$ belongs to $\Gamma_\R$.  
The uniqueness modulo $\F_l$ of $v_{l+1}$ implies that this notion is independent 
 of the choice of $v_l$ in Lemma \ref{l:F-l}. Observe that there are only finitely many almost $\Gamma$-radii.

\begin{lemma} \label{l:almost-Gamma}
Let $\Theta$ be a compact subset of $\R/2\pi\Z$ and let $\D_\Theta$ be the union of the radii $L_\theta$ with $\theta\in\Theta$. 
Assume that $\D_\Theta$ contains no almost $\Gamma$-radius. Then $\pi(f(x))$ has no cluster value in $\T$ when $x$ tends to $0$ and $x\in\D_\Theta$. 
In particular,  when there is no almost $\Gamma$-radius, or equivalently, when $\kappa>1$ and $\F_\kappa\cap \Gamma_\R=\F_{\kappa-1}$
or when $\kappa=1$ and $\F_\kappa\cap \Gamma_\R=0$ (see also the hypothesis (H2) below), we have $\Lambda_f=\varnothing$.
\end{lemma}
\proof
It is enough to show that $\dist(f(x),\Gamma_\R)$ tends to infinity when $x$ tends to 0 and $x\in\D_\Theta$. Since $v_1,\ldots,v_{\kappa-1}$ are in $\Gamma_\R$, using the expansion of $f$ in Lemma \ref{l:F-l}, we have 
$$\liminf_{x\to 0\atop x\in\D_\Theta} |x|^{d_\kappa} \dist(f(x),\Gamma_\R) =\inf_{\theta \in\Theta} \dist (e^{-id_\kappa\theta }v_\kappa,\Gamma_\R).$$
By hypothesis, $e^{-id_\kappa\theta}v_\kappa$ belongs to a compact set outside $\Gamma_\R$ when $\theta$ is in $\Theta$. Therefore, the last infimum is a strictly positive number.  The lemma follows.
\endproof

From now on, we assume the following hypothesis, see Lemma \ref{l:almost-Gamma}. 

\medskip\noindent
{\bf (H2)} There exists at least one almost $\Gamma$-radii or equivalently $\F_\kappa\cap \Gamma_\R$ is a real hyperplane of $\F_\kappa$. 

\medskip

Observe that the set of almost $\Gamma$-radii is the pullback of the real hyperplane $\F_\kappa\cap\Gamma_\R$ to $\D^*$ by the map $x\mapsto x^{-d_\kappa} v_\kappa$. This is a family of $2d_\kappa$ radii which are equidistributed on $\D$. Let $\lambda\in\C$ with $|\lambda|=1$ such that $\lambda^{d_\kappa} v_\kappa$ belongs to $\Gamma_\R$. The number $\lambda^{d_\kappa}$ does not depend on the choice of $v_l$ in Lemma \ref{l:F-l}. 

Choose a $\C$-linear function $H$ on $\E$ such that  the real hyperplane $\Im(H)=0$ contains $\Gamma_\R$ but not $\C v_\kappa$. This function is not unique in general. We have $\Im(H)=0$ on $\C v_1,\ldots,\C v_{\kappa-1}$  and hence $H=0$ on these complex lines. The restriction of $H$ to $\C v_\kappa$ is a non-zero linear polynomial and $\Im(H)=0$ on $\lambda^{d_\kappa}\R v_\kappa$. By multiplying $H$ by a real constant, we can assume that 
$H(\lambda^{d_\kappa} v_\kappa)=1$.  It follows that 
  $H(f(x))=(\lambda x)^{-d_\kappa}+O(x^{-d_\kappa+1})$. 
Then, there is a holomorphic function $x'=\lambda x+O(x)$, defined on some neighbourhood of 0 in $\D$, such that $H(f(x))=x'^{-d_\kappa}$. 

We will work with $x'$ as a new coordinate together with the above coordinate $x$. In what follows, the notation like $x_n'$ denotes the same point $x_n$ in the coordinate $x'$. Fix a number $\rho>0$ small enough such that the disc $\D'(\rho):=\{|x'|<\rho\}$ is contained in $\D$. 
Observe that the almost $\Gamma$-radii are $L_p:=\lambda^{-1} L_{p\pi/d_\kappa}$ with $p$ in $\{0,\ldots,2d_\kappa-1\}$. Define
$$L_p':=\Big\{x'=re^{{ip\pi\over d_\kappa}} \text{ with } r\in (0,\rho)\Big\}.$$
This is a radius of $\D'(\rho)$ which is tangent to $L_p$ at 0.

Fix  an arbitrary constant $A>0$ and consider the domain $\E_A$ in $\E$ defined by the inequality $|\Im(H)|<A$. 
We have the following lemma.

\begin{lemma} \label{l:D-H}
The set $f^{-1}(\E_A)\cap \D'(\rho)$ is the union of the following $2d_\kappa$ sets which are open near $0$
$$\Omega_{A,p}:=\Big\{x'=re^{{i(\theta+p\pi)\over d_\kappa}} \text{ with } r\in (0,\rho),\ |\theta|\leq {\pi\over 2} \text{ and } |\sin\theta| < A r^{d_\kappa}\Big\}$$
with $p=0,\ldots, 2d_\kappa-1$. In particular, if $x_n'$ is a sequence converging to $0$ such that $|x_n'|^{-d_\kappa-1}\dist(x_n',L'_p)$ tends to $\infty$ for every $p$, then the sequence $\pi(f(x_n'))$ has no cluster value in $\T$.
\end{lemma}
\proof
The first assertion is a direct consequence of the definitions of $x'$ and $\E_A$. When $x'\in\Omega_{A,p}$ with $r$ is small, we have $|\theta|\not=\pi/2$. So $\Omega_{A,p}$ is open near 0. 

For the second assertion, observe that $f(x_n')$ is outside $\E_A$ when $n$ is large enough. This is true for every $A$. It follows that the distance between $f(x_n')$ and $\Gamma_\R$ tends to infinity. Thus, the sequence $\pi(f(x_n'))$ has no cluster value in $\T$.
\endproof

\begin{lemma} \label{l:F-l-2}
Fix an integer $0\leq p\leq 2d_\kappa-1$. 
Then, there is a unique sequence of integers $\delta_0>\cdots>\delta_{m+1}$ with $\delta_0=d_\kappa, \delta_{m+1}=0$ and some $0\leq m\leq 2n-2\kappa+1$,  such that for $x'\in L_p'$  
$$f(x')=\sum_{j=1}^{\kappa-1} \big(1+ O(x')\big)\lambda^{d_j}x'^{-d_j}v_j+\sum_{j=0}^m \big(1+ O(x')\big) |x'|^{-\delta_j}\vartheta_j  +
\vartheta_{m+1} +O(x') \quad \text{as} \quad x'\to 0,$$
where 
$\vartheta_0,\ldots,\vartheta_{m+1}$ are vectors in $\E$ with $\vartheta_0,\ldots,\vartheta_m$ $\R$-linearly independent modulo $\F_{\kappa-1}$.
Denote by $\F_l'$ the real vector space of real dimension $2\kappa+l-1$, spanned by $\F_{\kappa-1}$ and $\vartheta_0,\ldots,\vartheta_l$. Then 
$\vartheta_0=\lambda^{d_\kappa} v_\kappa$ modulo $\F_{\kappa-1}$ and $\vartheta_{l+1}$ is unique modulo $\F_l'$. Moreover, we also have 
$\|f_{\F'_l}\|=\Theta(x'^{-{\delta_{l+1}}})$ for $l<m$ and $\|f_{\F_m'}\|=O(1)$, $\F'_m+i\F'_m= \F_k$ and  $V'\subset V$ with $V':=i\R \vartheta_0+ \F_m' + \vartheta_{m+1}$. 
\end{lemma}
\proof
Except for the last two inclusions, the lemma can be proved using the same arguments as in Lemma \ref{l:F-l}. 
We have $\delta_0=d_{\kappa}$ because both $x^{-d_{\kappa}}$ and $|x'|^{-\delta_0}$ appear in the leading term of $f_{\F_{\kappa-1}}$. 
As in the proof of  Lemma \ref{l:F-l}, $\F'_m$ is the minimal subspace of $\E$ such that $f_{\F'_m}$ is bounded on $L'_p$. 
So it is contained in $\F_k$ since  $f_{\F_k}$ is bounded on $\D^*$. Since $\F_k$ is complex, we deduce that $\F'_m+i\F'_m\subset \F_k$.
If we define $\F:=\F'_m+i\F'_m$, then $f_\F$ is bounded on $L'_p$. It follows that it is bounded on $\D^*$ and hence $\F\supset \F_k$. We conclude that $\F'_m+i\F'_m= \F_k$.

For the rest of the lemma, when $x'\in L'_p$ tends to 0, the distances from $f(x')$ to both $\F_m' + \vartheta_{m+1}$ and $V$ tend to 0. 
Since $\F_m' + \vartheta_{m+1}$ and $V$ are obtained from $\F_m'$ and $\F_k$ by translation, we deduce that 
 $\F_m' + \vartheta_{m+1}$ is contained in $V$. Finally, as $v_\kappa$ is contained in $V$, the complex line $\C\vartheta_0$ is also contained in $V$ and hence $V'$ is contained in $V$.
\endproof

Note that $\delta_l, \vartheta_l$ and $\F_l'$ depends on $L_p$. 
We say that  $L_p'$ is {\it a $\Gamma$-radius} if $\F_m'$ is contained in $\Gamma_\R$. 
In this case, since $\vartheta_{l+1}$ is unique modulo $\F_l'$, the vector $\vartheta_l$ belongs to $\Gamma_\R$ for every $l\leq m$. 

\begin{lemma} \label{l:Gamma}
Let $\Theta$ be a compact subset of $\R/2\pi\Z$ and let $\D'_\Theta$ be the union of the radii $L'_\theta$ with $\theta\in\Theta$. 
Assume that $\D'_\Theta$ contains no $\Gamma$-radius. Then $\pi(f(x'))$ has no cluster value when $x'\in\D'_\Theta$ and $x'$ tends to $0$.
In particular,  when there is no  $\Gamma$-radius (see also the hypothesis (H3) below), we have $\Lambda_f=\varnothing$.
\end{lemma}
\proof
By Lemma \ref{l:almost-Gamma}, it is enough to consider $\Theta$ such that $\D_\Theta$ is a small sector containing one and only one almost $\Gamma$-radius $L_p$ such that $L_p'$  is not a $\Gamma$-radius. Let $x'_j$ be a sequence in $\D'(\rho)$ converging to $0$ in the direction of $L'_p$ such that $|x'_j|^{-d_\kappa-1}\dist(x'_j,L'_p)$ is bounded by a constant. By Lemma \ref{l:D-H}, we only have to check that $\dist (f(x'_j),\Gamma_\R)$ tends to infinity. 
Let $l\geq 1$ be the smallest integer such that $\vartheta_l$ does not belong to $\Gamma_\R$. 

Let $y_j'$ be the projection of $x'_j$ to $L'_p$. We have $|x'_j-y_j'| = O(x_j'^{d_\kappa+1})$ and hence 
$x_j'=y_j'(1+O(y_j'^{d_\kappa}))$.
Using the expansion of $f(x)$ in Lemma \ref{l:F-l} and the fact that $v_1,\ldots,v_{\kappa-1}$ belong to $\Gamma_\R$, we obtain 
$$\dist(f_{\F_{\kappa-1}}(x_j'),f_{\F_{\kappa-1}}(y_j'))=O(1)  \quad \text{and hence}\quad \dist(f(x_j'),\Gamma_\R)=\dist(f(y_j'),\Gamma_\R)+O(1).$$
Here, we use that the factors $1+O(x)$ in  Lemma \ref{l:F-l} are analytic functions in $x$.

Using the expansion of $f(x)$ in Lemma \ref{l:F-l-2}, applied to $y_j'$ instead of $x'$, and the fact that
 $v_1,\ldots, v_{\kappa-1},\vartheta_0,\ldots, \vartheta_{l-1}$ belong to $\Gamma_\R$, we have 
$$|y_j'|^{\delta_l}\dist(f(y_j'),\Gamma_\R) = \dist(\vartheta_l, \Gamma_\R) + O(y_j').$$
When $x_j'$ is small enough, the last expression is positive since $\vartheta_l$ does not belong to $\Gamma_\R$. 
Therefore, $\dist(f(y_j'),\Gamma_\R)$ tends to infinity.
The lemma follows.
\endproof

From now on, we assume the following hypothesis, see Lemma \ref{l:Gamma}.
 
\medskip\noindent
{\bf (H3)}  The exists  at least one $\Gamma$-radius.

\medskip

We will use the notations introduced in Lemma \ref{l:F-l-2}. 
Let $\Theta$ be a compact subset of $\R/2\pi\Z$. Assume that $\D_\Theta'$ contains only one  $\Gamma$-radius $L'_p$, and this radius is in the interior of $\D'_\Theta$. Denote by
$\Lambda_{f,\Theta}$ the set of all cluster values of $\pi(f(x))$ for $x\in \D_\Theta$ going to 0. By Lemma \ref{l:almost-Gamma}, this is also  the set of all cluster values of $\pi(f(x'))$ for $x'\in \D_\Theta'$ going to 0. 

\begin{lemma} \label{l:V-prime}
We have that $\Lambda_{f,\Theta}$ is contained in $\overline{\pi(V')}$. 
\end{lemma}
\proof
Let $x_j$ be a sequence in $\D_\Theta$ converging to 0 such that $\pi(f(x_j))$ converges to some point in $\T$. We deduce that the distance between $f(x_j)$ and $\Gamma_\R$ is bounded by a constant independent of $j$. 
We will use as above the notation $x_j'$ for the same point $x_j$ in coordinate $x'$.
By Lemma \ref{l:D-H}, we have $\dist(x_j',L'_p)=O(x_j'^{d_\kappa+1})$.

We use  the notation from the proof of Lemma \ref{l:Gamma}. Since both $\F_\kappa+\vartheta_{m+1}$ and  $\F_m'+\vartheta_{m+1}$ are contained in $V'$, we get 
$$\dist (f(x'_j), V')=\dist (f(y_j'), V') +O(x_j')=O(x_j').$$
Therefore, the limit of $\pi(f(x'_j))$ belongs to  $\overline{\pi(V')}$. The lemma follows. 
\endproof

The lemma below gives us a characterization of sequences of values of $\pi\circ f$ having  cluster points in $\T$, see also Lemma \ref{l:D-H}.

\begin{lemma} \label{l:Omega}
 Let $x_j$ be a sequence of points  in $\D_\Theta$ going to $0$ and let $x_j'$ be the same sequence written in coordinate $x'$. Then the sequence $\pi(f(x_j'))$ is relatively compact in $\T$ if and only if $x_j'^{-d_\kappa-1}\dist(x_j',L'_p)$ is a bounded sequence. In particular, the image of $\Omega_{A,p}$ by $\pi\circ f$ is contained in a compact subset of $\T$.
\end{lemma}
\proof
Assume that $x_j'^{-d_\kappa-1}\dist(x_j',L'_p)$ is a bounded sequence. 
By Lemma \ref{l:D-H}, it is enough to show that $\dist(f(x_j'),\Gamma_\R)$ is bounded. Arguing as in the proof of Lemma \ref{l:Gamma}, we only need to check that 
$\dist(f(y_j'),\Gamma_\R)$ is bounded. This is clear because $f(y_j')$ tends to $\F_m'+\vartheta_{m+1}$ which is parallel to $\Gamma_\R$ as $L_p'$ is a $\Gamma$-radius.
The lemma follows.
\endproof

We can now study $\Lambda_{f,\Theta}$ as we did for $\Lambda_f$ in Section \ref{s:compact-1D}. 
Let $\U$ be a relatively compact open subset of $\D^*$ with piecewise smooth boundary. 
Assume that $\U$ is contained in $\Theta$ and its boundary intersects $L_p$ transversally. Assume also that $\U\cap L_p$ is non-empty, see Lemma \ref{l:almost-Gamma}.
For $a\in (0,1] $ define $\U_a:=a\U$.  Notice that unlike the previous section, here $a$ is a real positive number. 
The following lemma is similar to Lemma \ref{l:mass}. 

\begin{lemma} \label{l:mass-2}
There is a constant $\lambda_A>0$ such that the mass in $\E_A$  of the positive measure $f_* [\U_a]\wedge \omega_\E$ is equal to $\lambda_A a^{-2d+d_\kappa} + O(a^{-2d+d_\kappa+1})$ when $a$ tends to $0$.
\end{lemma}
\proof  
Let $\U_{A,a}$ denote the intersection of $\U_a$ with $\Omega_{A,p}$ and define $\widetilde \U_{A,a}:=a^{-1}\U_{A,a}$. As in Lemma \ref{l:mass}, the mass of $f_* [\U_a]\wedge \omega_\E$ in $\E_A$ satisfies
\begin{eqnarray*}
\|f_* [\U_a]\wedge \omega_\E \|_{\E_A}  & = &    d^2 \|v_1\|^2  \int_{x\in \U_{A,a}} \big(|x|^{-2d-2}+O(x^{-2d-1})\big) (i dx \wedge d \overline x)\\
& = & d^2 \|v_1\|^2  a^{-2d} (1+O(a)) \int_{x\in \widetilde \U_{A,a}} |x|^{-2d-2} (i dx \wedge d \overline x).
\end{eqnarray*}

Now, the intersection of $\U$ with $L_p$ is a finite union of intervals that we denote by $(q_1,q_1'),\ldots,(q_s,q_s')$. The tangent cones of the boundary of $\U$ at 
$q_1,q_1'\ldots,q_s,q_s'$ are denoted by $\gamma_1,\gamma_1'\ldots,\gamma_s,\gamma_s'$. They are real lines or unions of two half-lines. We see that $\widetilde \U_{A,a}$ is approximatively the union of $p$ open subsets of $a^{-1}\Omega_{A,p}$ limited by $\gamma_j$ and $\gamma_j'$. By considering the boundary of $a^{-1}\Omega_{A,p}$, we see that the last integral is equal to
$$\beta_A a^{d_\kappa} +O(a^{d_\kappa+1}),$$
for some constant $\beta_A>0$. The lemma follows. Note that the length of the part of the boundary of $\U_{A,a}$ inside $\Omega_{A,p}$ (resp. of $\widetilde \U_{A,a}$ inside $a^{-1}\Omega_{A,p}$) is $O(a^{d_\kappa+1})$ (resp. $O(a^{d_\kappa})$) and $\beta_A$ depends continuously on $A$. We will use these properties later.
\endproof

Define
$$\mu_a:= a^{2d-d_\kappa} \pi_*\big(f_* [\U_a] \wedge \omega_\E\big)$$
which is a positive measure on $\T$. 

\begin{theorem} \label{t:1D-Gamma}  
We use the above notations and  assume that the hypotheses (H0)-(H3) hold (otherwise, $\Lambda_f$ is empty, a singleton, or Theorem \ref{t:compact-1D} applies). Then, the measure $\mu_a$ converges to a Haar measure of $\overline{\pi(V')}$ as $a$ tends to $0$. Moreover, we have 
$\Lambda_{f,\Theta}=\overline{\pi(V')}$ and $\Lambda_f$ is a finite union of translated connected closed real Lie subgroups of $\T$. 
\end{theorem}

Note that the second assertion is a direct consequence of the first one and Lemma \ref{l:V-prime}. So we only need to prove the convergence of $\mu_a$ in the last theorem. By Lemma \ref{l:F-l-2}, we also notice that $V'$ is obtained from the vector space $i\R\vartheta_0+\F_m'$ by a translation and $i\R\vartheta_0+\F_m'$ is contained in the complex vector space spanned by its intersection with $\Gamma_\R.$ We will use this property in the last section for the proof of Theorem \ref{t:ALW}.

As in Section \ref{s:compact-1D}, for $1\leq l\leq \kappa$, 
we define the map $\tau_{l,a}:\overline \U\times \overline \D^{l-1}\to \E$ by
$$\tau_{l,a}(x,t_1,\ldots,t_{l-1}) := f(ax)+ \sum_{j=1}^{l-1} a^{-d_j}t_j v_j$$
and the positive current 
$S_{l,a}$ of bi-dimension $(l,l)$ on $\E$ by 
$$S_{l,a}:= c_{l,a}^{-1}(\tau_{l,a})_*\big(\Psi_l [\U\times\D^{l-1}]\big),$$
where 
$$\Psi_l(x,t_1,\ldots,t_{l-1}):=|x|^{-2(d-d_l)} \qquad \text{and} \qquad c_{l,a}:=  (2\pi)^{l-1}d_l^2a^{-2d_1-\cdots-2d_l+d_\kappa}.$$
Note that the definitions of the current $S_{l,a}$ and the constant $c_{l,a}$ here are slightly different from the ones in Section \ref{s:compact-1D}, due to an extra factor $a^{d_\kappa}$ in $c_{l,a}$.

Recall that $\I:=(0,1)$. For $1\leq l\leq m$, define the map $\tau_{l,a}':\overline\U\times\overline\D^{\kappa-1}\times \overline \I^l\to \E$ by 
$$\tau'_{l,a}(x,t_1,\ldots,t_{\kappa-1},s_0,\ldots,s_{l-1}):= f(ax)+ \sum_{j=1}^{\kappa-1} a^{-d_j}t_j v_j + \sum_{j=0}^{l-1} a^{-\delta_j}s_j \vartheta_j$$
and the current $S_{l,a}'$ of dimension $2\kappa+l$ on $\E$ by 
$$S_{l,a}':=(c_{l,a}')^{-1}(\tau_{l,a}')_*\big(\Psi_l' [\U\times\D^{\kappa-1}\times \I^l]\big),$$
where 
$$\Psi_l'(x,t_1,\ldots,t_{\kappa-1},s_0,\ldots,s_{l-1}):=|x|^{-(2d-d_\kappa-\delta_l)}$$
and
$$c_{l,a}':=  (2\pi)^{\kappa-1}d_\kappa\delta_l a^{-2d_1-\cdots-2d_{\kappa-1}-d_\kappa-\delta_1-\cdots-\delta_l}.$$

For simplicity, define $t:=(t_1,\ldots,t_{\kappa-1})$ and $s:=(s_0,\ldots,s_{l-1})$.
The following lemma will be used together with Lemma \ref{l:vector}.

\begin{lemma} \label{l:vector-2}
Let  $\eta_0,\ldots,\eta_{l-1}$ be the images of the tangent vectors
$${\partial\over\partial s_0} \CommaBin \   \cdots \  \CommaBin \  {\partial\over\partial s_{l-1}} $$
by the differential $d\tau'_{l,a}(x,t,s)$ of $\tau'_{l,a}$ at the point $(x,t,s)$. 
Define 
$$\sigma_l:=2d_1+\cdots+2d_\kappa+\delta_1+\cdots+\delta_{l-1}.$$ 
Then we have for $1\leq l\leq m$ and $ax\in \U_a\cap\Omega_{A,p}$
\begin{multline*}
(iu\wedge \overline u) \wedge (iu_1 \wedge \overline u_1) \wedge \ldots \wedge  (iu_{\kappa-1} \wedge \overline u_{\kappa-1}) \wedge \eta_0 \wedge \ldots  
\wedge \eta_{l-1} = \\
d_\kappa\delta_l a^{-\sigma_l}  |x|^{-d_\kappa-\delta_l-2} (iv_1\wedge\overline v_1)\wedge\ldots\wedge (iv_{\kappa}\wedge \overline v_{\kappa})\wedge \vartheta_1\wedge\ldots\wedge \vartheta_l  +O(a^{-\sigma_l+1}).
\end{multline*}
Moreover, if $w_1,\ldots,w_{2\kappa+l-1}$ belong to the family of vectors $u,\overline u,u_1,\overline u_1,\ldots,\eta_{l-1}$, then
$$w_1\wedge \ldots\wedge w_{2\kappa+l-1} = O(a^{-\sigma_l+1}).$$
\end{lemma}
\proof
Note that by Lemmas \ref{l:F-l} and \ref{l:F-l-2}
$$(iv_1\wedge\overline v_1)\wedge\ldots\wedge (iv_{\kappa}\wedge \overline v_{\kappa})=
(iv_1\wedge\overline v_1)\wedge\ldots\wedge (iv_{\kappa-1}\wedge \overline v_{\kappa-1}) (i\vartheta_0\wedge\overline\vartheta_0).$$ 
Denote by $u'$ the image of $\partial/\partial x'$ by $d\tau'_{l,a}(x,t,s)$. Then we have $iu'\wedge\overline u'= (1+O(a))iu\wedge \overline u$. So we can replace $u,\overline u$ by $u',\overline u'$.  
Recall that the points in $L_p'$ have arguments $p\pi/d_\kappa$ with respect to the coordinate $x'$. 
Therefore,  by Lemma \ref{l:F-l-2} and the uniqueness theorem for holomorphic functions, we have for $x'\in \D'(0,\rho)$
\begin{align} \label{e:xprime-f-vartheta}
f(x')=\sum_{j=1}^{\kappa-1} \big(1+ O(x')\big)\lambda^{d_j}x'^{-d_j}v_j+\sum_{j=0}^m \big(1+ O(x')\big) e^{i\delta_j p\pi/d_\kappa} x'^{-\delta_j}\vartheta_j  + \vartheta_{m+1} + O(x').
\end{align}

Recall that for $ax\in \U_a\cap\Omega_{A,p}$ we have $\arg(x)=p\pi/d_k + O(a)$. Then,
using the formula \eqref{e:xprime-f-vartheta}, we obtain
$$u'=-\delta_l a^{-\delta_l} e^{i p \pi/d_\kappa} |x'|^{-\delta_l-1}\vartheta_l+O(a^{-\delta_l+1}) \quad \text{modulo} \quad v_1,\overline v_1,\ldots, v_{\kappa-1}, \overline v_{\kappa-1}, \vartheta_0,\ldots,\vartheta_{l-1}$$ 
and 
$$\overline u'= -d_\kappa a^{-d_\kappa}  e^{-i p\pi/d_\kappa}|x'|^{-d_\kappa-1}\overline\vartheta_0 + O(a^{-d_\kappa+1})  \quad \text{modulo} \quad  
v_1,\overline v_1,\ldots, v_{\kappa-1}, \overline v_{\kappa-1}.$$ 
Now, we just need to follow the proof of Lemma \ref{l:vector} in order to get the result. 
\endproof

The following lemma is a version of Lemma \ref{l:limit-value-S-l}.

\begin{lemma} \label{l:limit-value-S-l-2}  
All limit values of $\pi_*(S_{l,a})$ (resp. $\pi_*(S'_{l,a})$), when $a$ tends to $0$, are non-zero positive (resp. non-zero) closed currents of bi-dimension $(l,l)$ (resp. dimension $2\kappa+l$)
invariant by $\F_l$ (resp. $\F_l'$). Moreover, they have the same mass on $\T_A:=\pi(\E_A)$ for every $A$.
\end{lemma}
\proof 
The proof follows the one of Lemma \ref{l:limit-value-S-l} using Lemmas \ref{l:vector} and \ref{l:vector-2}. We will only explain the difference in the new setting.
Consider the case of $\pi_*(S_{l,a})$. The proof for $\pi_*(S'_{l,a})$ is similar.
Since $\E_A$ is  a union of affine spaces parallel to $\Gamma_\R$, we have $\E_A=\pi^{-1}(\T_A)$.  By the last assertion of Lemma \ref{l:Omega},  the intersections of the supports of  $S_{l,a}, S'_{l,a}$ with $\E_A$ are sent by $\pi$ to a compact subset of $\overline\T_A$ which does not depend on $a$.

Let $\widehat S_{l,a}$ be the restriction of $S_{l,a}$ to $\E_A$. The restriction of $\pi_*(S_{l,a})$ to $\T_A$ is equal to $\pi_*(\widehat S_{l,a})$.
In this lemma, we only need to consider currents and their boundaries inside $\E_A$ and $\T_A$. 

We show that the mass of $\widehat S_{l,a}$ is equal to $\beta_{l,A}+O(a)$ for some constant $\beta_{l,A}$ depending continuously on $A$. 
The continuity of $\beta_{l,A}$ in $A$ insures that the limit currents of  $\pi_*(\widehat S_{l,a})$ have no mass on the boundary of $\T_A$ and hence have mass $\beta_{l,A}$ on $\T_A$.
We also show that the mass of $d\widehat S_{l,a}$ in $\E_A$ is $O(a)$. Then, the result will follow  as in the proof of Lemma \ref{l:limit-value-S-l}.
We use here the notations from the proof of Lemma \ref{l:mass-2}. Observe that the support of $\widehat S_{l,a}$ is the image of $(a^{-1}\overline\U_{A,a})\times \overline \D^{l-1}$ by $\tau_{l,a}$. Here, $\overline \U_{A,a}$ is the closure of $\U_{A,a}$ in $\Omega_{A,p}$. We also need to consider the boundary of $\U_{A,a}$ but only the part inside $\Omega_{A,p}$.

The volume of $(a^{-1}\overline\U_{A,a})\times \overline \D^{l-1}$ is $O(a^{d_\kappa})$ and the length of the part of  boundary of $a^{-1}\overline\U_{A,a}$ that we consider is also $O(a^{d_\kappa})$. This is the main difference in comparison with the situation in Lemma \ref{l:limit-value-S-l}. This is why in the definition of $S_{l,a}$ in this section, the  constant  $c_{l,a}$ is  different from the one in Section \ref{s:compact-1D} by a factor $a^{d_\kappa}$. Now, we can obtain the result using the arguments from Lemma \ref{l:limit-value-S-l}. As in Lemma \ref{l:mass-2}, it is not difficult to see that the involving constant $\beta_{l,A}$ depends continuously on $A$. 
\endproof

The following lemma is obtained in the same way as for Proposition \ref{p:limit-S-k}  using Lemmas \ref{l:vector}, \ref{l:vector-2} and \ref{l:limit-value-S-l-2}.

\begin{proposition} \label{p:limit-S-m-1}  
When $a$ tends to $0$, the current $\pi_*(S'_{m,a})$ tends to a Haar current  associated with $V'$ that we will denote by 
$R'_{m}$.
\end{proposition}

Choose $m$ differential 1-forms $\Psi_1,\ldots,\Psi_{m}$ with constant coefficients on $\E$ such that $\Psi_l(v_j)=0$,  $\Psi_l(\overline v_j)=0$ for $1\leq j\leq \kappa$, and $\Psi_l(\vartheta_j)=1$ if $j=l$ and $\Psi_l(\vartheta_j)=0$ if $j\not=l$.
The following proposition is obtained in the same way as for Proposition \ref{p:limit-S-l}.

\begin{proposition} \label{p:limit-S-l-2}
For every $1\leq l\leq m$, the current
$\pi_*(S'_{l, a})$ tends to the closed current 
$$R'_l:= R'_{m}  \wedge  \Psi_{l+1} \wedge \ldots \wedge \Psi_{m},$$ 
as $a$ tends to $0$. For every $1\leq l\leq \kappa$, the current
$\pi_*(S_{l, a})$ tends to the positive closed current 
$$R_l:= R'_{m}  \wedge (i \Phi_{l+1} \wedge \overline \Phi_{l+1})\wedge \ldots \wedge (i \Phi_{\kappa} \wedge \overline \Phi_{\kappa}) \wedge  \Psi_1 \wedge \ldots \wedge \Psi_{m}$$ 
as $a$ tends to $0$. Moreover, the supports of $R_l$ and $R_l'$ are both equal to $\overline{\pi(V')}.$  
\end{proposition}

\proof[End of the proof of Theorem \ref{t:1D-Gamma}]  It is enough to follow the arguments at the end of the proof of Theorem \ref{t:compact-1D}
and obtain that $\mu_a$ converges to a Haar measure on 
 $\overline{\pi(V')}$. As mentioned just after Theorem \ref{t:1D-Gamma}, this convergence property implies  the second assertion in this theorem.
\endproof

Note that in general $\Lambda_f$ may not be irreducible. Consider the following example. Let $\E:=\C^2$, $\Gamma:=\Z\times (\Z+i\Z)$ and 
$\T:=\E/\Gamma=(\C/\Z)\times (\C/(\Z+i\Z))$. Consider $f:\overline\D^*\to\E$ defined by $f(x):=(x^{-2},x^{-1})$. 
The almost $\Gamma$-radii are $L_p:=L_{p\pi/2}$ with $p=0,1,2,3$.
We can take $x'=x$ and see that all $L_p$ are $\Gamma$-radii. For $p=0,2$, we have $V'=\C\times \R$ and for $p=1,3$ we have $V'=\C\times (i\R)$. 
In this example, we get two different real semi-tori $\overline{\pi(V')}=\C^*\times (\R/\Z)$ when $p=0,2$ and $\overline{\pi(V')}=\C^*\times (i\R/i\Z)$ when $p=1,3$.

\proof[Proof of Theorem \ref{t:main-1}]  
We identify the vector space $\E$  to an affine chart of the projective space $\overline\E:=\P(\E\oplus \C)$ in the natural way. 
Observe that the usual  closure $\overline X$ of $X$ in this projective space is an algebraic curve which is the union of $X$ with finitely many  points at the hyperplane at infinity $\overline\E\backslash\E$. 
Any irreducible germ of $\overline X$ at a point $x_\infty$ at infinity can be parametrized by a smooth map $f:\overline\D\to \overline\E$ such that 
$f$ is holomorphic in $\D$, $f(\overline\D^*)\subset X$ and $f(0)=x_\infty$. This map has a polynomial growth in the sense of Section \ref{s:compact-1D} when $x$ tends to 0. Now, it is enough to apply Theorem \ref{t:1D-Gamma} and get the result.
\endproof

\section{Higher dimensional algebraic flows in a torus} \label{s:higherdim}

In this section, we will prove a slightly more general version of Theorem \ref{t:main-2} for flows in a compact torus, see Theorem \ref{t:main_high_compact} below. The result still holds when $\T$ is non-compact provided that the flow stays inside a compact set. At the end of this section, we will consider an example of a flow with non-compact support.

\medskip\noindent
{\bf Main result and preparation step of the proof.} 
Let $\Omega$ be an irreducible complex variety of dimension $q$ and let $Y$ be a proper compact analytic subset of $\Omega$. Consider a meromorphic map 
$f:\Omega\to \E$ which is holomorphic on $\Omega\backslash Y$. That is, $f(z)$ has a polynomial growth, with respect to the inverse of the distance from $z$ to $Y$, when $z$ tends to $Y$.  As before, our goal is to describe the set $\Lambda_f$  of limit points of $\pi\big(f(z)\big)$ in $\T:=\E/\Gamma$ as $z \in \Omega \backslash Y$ tends to $Y$.
 The main result of this section is the following.

\begin{theorem} \label{t:main_high_compact}  
Let $\Omega,Y$ and $f$ be as above and assume moreover that $\T=\E/\Gamma$ is compact. 
Let $Y^*$ be an analytic subset of $Y$ and $\Lambda_{f,Y^*}$ the set of limit points of $\pi\big(f(z)\big)$ in $\T$ as $z \in \Omega \backslash Y$ tends to $Y^*$.
Then, 
there exist finitely many (possibly trivial) complex vector subspaces $\F_1,\ldots,\F_m$ of $\E$  and  complex algebraic subsets $C_1,\ldots,C_m$ of $\E$ of dimension less than $\dim \Omega$, all independent of $\Gamma$,  such that
\begin{align} \label{e:Lambda-decomp}
\Lambda_{f,Y^*}= \bigcup_{j=1}^m \big(\pi(C_j)+ \T_j \big).
\end{align}
Here, $\T_j$ denotes the closure of $\pi(\F_j)$ in $\T$ and is a real torus.
Moreover, if $\F_j$ is maximal among $\F_1,\ldots,\F_m$ for the inclusion,  then $C_j$ is a finite set. 
\end{theorem}

Note that Theorem \ref{t:main-2} is a direct consequence of the last theorem. Indeed, we can naturally identify $\E$ with an affine chart of the projective space $\P(\E \oplus \C)$. Then, we just need to take $\Omega$ to be the closure of $X$ in $\P(\E \oplus \C)$, $Y^*:=Y:=\Omega\backslash X$ 
and $f:\Omega\backslash Y\to \E$  the natural inclusion map.  Note that in this setting, we will see during the proof below that the tori $\T_i$ are non-trivial as $f(z)$ tends always to infinity when $z$ tends to $Y$.

Observe that if $\tau:\widetilde\Omega\to\Omega$ is a composition of blow-ups with smooth centers, $\widetilde Y:=\tau^{-1}(Y)$ and $\widetilde Y^*:=\tau^{-1}(Y^*)$, then we can replace $\Omega,Y,Y^*$ and $f$ by $\widetilde\Omega,\widetilde Y,\widetilde Y^*$ and $f\circ\tau$ because  $\Lambda_f=\Lambda_{f\circ\tau}$. Therefore, using a resolution of singularities and other suitable blow-ups, we can assume, from now on, that $\Omega$ is smooth, $Y$ is a finite union of compact smooth hypersurfaces having simple normal crossings everywhere and $Y^*$ is a union of some irreducible components of $Y$, see \cite{Hironaka}. 

Let $\Rf=\{Y_1, \ldots, Y_l\}$  be a family of irreducible components of $Y$ such that $Y_\Rf:=\cap_{j=1}^l Y_j$ is non-empty. Denote by
$Y_\Rf^c$ the union of the irreducible components of $Y$ which are not contained in $\Rf$.  
Let $Z$ be any irreducible component of $Y_\Rf$. This is a smooth compact submanifold of codimension $l$ of $\Omega$.
 Let $\Lambda_{f,Z}$  be the set of limit points of $\pi\big(f(z)\big)$ when $z$ tends to $Z\backslash Y_\Rf^c$. Observe that $\Lambda_{f,Y^*}$ is a finite  union of  considered sets $\Lambda_{f,Z}$. 
So we only need to prove that $\Lambda_{f,Z}$ admits a decomposition similar to the right-hand side of \eqref{e:Lambda-decomp}. We need to introduce some auxiliary notations and results.

Consider a small non-empty open subset $U$ of $Z$ and a local holomorphic system of coordinates $(z_1, \ldots, z_q)$ in a neighbourhood of $\overline U$ in $\Omega$  on which 

\medskip\noindent
{\bf (S1)} the hypersurface $Y_j$ is given by the equation $z_j=0$ for $1\leq j\leq l$;

\medskip\noindent
{\bf (S2)} $Y^c_\Rf$ is a union, possibly empty, of the fibres of the map $(z_1,\ldots,z_q)\mapsto (z_{l+1},\ldots,z_q)$.

\medskip

Put $z':= (z_1, \ldots, z_l)$ and $z'':= (z_{l+1}, \ldots, z_q)$. We can use $z''$ as a coordinate system for $U$. For every $l$-tuple $\beta=(\beta_1, \ldots, \beta_l) \in \Z^l$, define 
$$|\beta|:= |\beta_1|+ \cdots+ |\beta_l|, \quad z'^\beta:= z^{\beta_1}_1 \ldots z^{\beta_l}_l
 \quad \text{and} \quad |z'|^\beta:= |z_1|^{\beta_1} \ldots |z_l|^{\beta_l}$$
(note that the last formula is also meaningful for $\beta\in \R^l$).  
For $\beta, \beta' \in \R^l,$ we write $\beta \ge \beta'$ and $\beta' \le \beta$  if  every component of $\beta -\beta'$ is nonnegative.  
For such $\beta,\beta'$, we also write $\beta > \beta'$ and $\beta'  < \beta$ if moreover $\beta\not=\beta'$.

Since $f$ is meromorphic with poles in $Y$, the Laurent expansion of $f(z)$ in $z'$ has the form
\begin{equation} \label{e:Laurent}
f(z)=\sum_{\beta\in\Z^l} z'^\beta v_\beta,
\end{equation}
where $v_\beta$ is a meromorphic function on $U$ with values in $\E$ which is holomorphic outside the sets $U\cap Y^c_\Rf$.
Observe  that if $v_\beta\not=0$, then the components of $\beta$ are bounded from below by a fixed negative constant which is determined by the order of the poles of $f$. 

\medskip\noindent
{\bf Leading powers and the vector space $\F_\beta$.} The cluster values of $\pi(f(z))$ depend on the "direction" that $z$ tends to $Z\backslash Y^c_\Rf$. In order to understand the different behaviors of $\pi(f(z))$ we will use a notion of complete leading sequences of powers. We first introduce the notion of leading powers.

For any power $\beta$, denote by $\F_\beta$ the complex vector subspace of $\E$ spanned by the following uncountable family of vectors
$$v_\beta(z'') \quad \text{with}\quad z''\in U\backslash Y^c_\Rf.$$
A power $\beta$ is said to be {\it leading} for $f$ if 

\medskip\noindent
{\bf (L1)} $v_\beta\not= 0$, at least one component of $\beta$ is negative, and 
$v_{\beta'}=0$ for every $\beta'<\beta$.

\begin{lemma} \label{l:dom_power}
We have the following properties.
\begin{itemize}
\item[(1)] There are finitely many leading powers for $f$.
\item[(2)] There is at least one leading power for $f$ unless $f$ is locally bounded near each point of $Z\backslash Y^c_\Rf$.
\item[(3)]  The set of leading powers of $f$ doesn't depend on $U$ nor on the chosen coordinate system (S1)-(S2).
\item[(4)] If $\beta$ is a leading power for $f$ as above, then $\F_\beta$ doesn't depend on $U$ nor on the chosen coordinate system (S1)-(S2).
\end{itemize}
\end{lemma}
\proof
(1) We only consider $\beta$ such that $v_\beta\not=0$. Since $f$ is meromorphic, the coefficients of such a $\beta$ are bounded from below by a fixed constant. Consider now only the minimal powers $\beta$ among those with $v_\beta\not=0$. It is enough to show that this family is finite. 

Assume by contradiction that this family contains an infinite sequence $\beta^{(1)},\beta^{(2)},\beta^{(3)}\ldots$ of distinct powers such that no one is strictly smaller than another. In particular, we don't have any subsequence whose coefficients tend to infinity.  So after extracting a subsequence, we can assume that one of the coefficients of $\beta^{(j)}$ does not depend on $j$.
Therefore, we can remove this coefficients from $\beta^{(j)}$, repeat the above argument a finite number of times and reach a contradiction.

\smallskip

(2) Assume that $f$ is not locally bounded near each point of $Z\backslash Y^c_\Rf$. So $Z$ is contained in the pole set of $f$. It follows that $v_{\beta^{(0)}}\not=0$ for some power $\beta^{(0)}$ which has at least a negative coefficient. 
Recall that when $v_\beta\not=0$, the coefficients of $\beta$ are bounded from below by a fixed negative number. Therefore, we can
choose a $\beta$ minimal with $v_\beta\not=0$ and $\beta\leq \beta^{(0)}$. Clearly, this $\beta$ is a leading power.

\smallskip

(3) Let $\widetilde U$ be another open subset of $Z$ and let $\widetilde z=(\widetilde z',\widetilde z'')$ be another coordinate system in a neighbourhood of $\widetilde U$ 
satisfying the conditions (S1)-(S2) above. We first consider the case where $W:=U\cap \widetilde U$ is non-empty. So, in a neighbourhood of $W$, $z$ and $\widetilde z$ denote the same point in two different coordinate systems. We have the following Laurent expansion
\begin{equation} \label{e:Laurent-bis}
f(\widetilde z)=\sum_{\beta\in\Z^l} \widetilde z'^\beta \widetilde v_\beta,
\end{equation}
where $\widetilde v_\beta$ is a meromorphic function on $\widetilde U$ with values in $\E$ which is holomorphic outside the set $\widetilde U\cap Y^c_\Rf$.

Observe that for $1\leq j\leq l$, since the equations $z_j=0$ and $\widetilde z_j=0$ define the same hypersurface, we have 
$$z_j= h_j \widetilde z_j +O(\widetilde z_j^2),$$
where $h_j$ is a nowhere vanishing holomorphic function on $W$. Define $h':=(h_1,\ldots,h_l)$. The relations between $z_j$ and $\widetilde z_j$, together with \eqref{e:Laurent}, give us
\begin{equation}  \label{e:Laurent-ter}
f(\widetilde z)=\sum_{\beta\in\Z^l} \widetilde z'^\beta h'^\beta  v_\beta + \text{higher order terms.}
\end{equation}
Here, by higher order terms, we mean terms involving powers strictly larger than  some $\beta$ with $v_\beta\not=0$. 

Clearly, the so-called higher order terms here cannot be leading and we see that the leading powers are the same for both coordinate systems $z$ and $\widetilde z$. Note that the last property still holds when $U\cap\widetilde U=\varnothing$ because we can connect $U$ and $\widetilde U$ using a chain of small open subsets in $Z$ and apply the previous case to consecutive open sets in the chain.

\smallskip

(4) Consider the same situation as above with $W\not=\varnothing$ and assume that $\beta$ is a leading power. 
First, observe that by uniqueness theorem for holomorphic functions, in the definition of $\F_\beta$, we can replace $U$ by $W$ and still get the same vector space.

Now, using 
\eqref{e:Laurent-bis} and \eqref{e:Laurent-ter}, we get 
\begin{equation} \label{e:coord-change}
\widetilde v_\beta =  h'^\beta  v_\beta \quad \text{on} \quad W.
\end{equation}
It follows that the vectors $v_\beta(a)$ and $\widetilde v_\beta(a)$ are co-linear for each point $a\in W\backslash Y^c_\Rf$. 
Therefore, the space $\F_\beta$ defined using $W, z$ coincides with the one defined by using $W,\widetilde z$.
The result follows.
\endproof

\noindent
{\bf Leading sequences of powers, $\F_\Bf$, $\E_\Bf$, $f_\Bf$, $\Sigma^-_\Bf$ and $\Sigma^{\geq 0}_\Bf$.} 
Consider a (possibly empty) finite sequence of powers $\Bf=(\beta^{(1)},\ldots, \beta^{(m)})$ with $\beta^{(j)} \in \Z^l$. Define 
$\F^{(0)}_\Bf:=0$, and for $1\leq j\leq m$, define the spaces $\F_\Bf^{(j)}$, $\F_\Bf$, $\E_\Bf$ and the map $f_\Bf:\E\to\E_\Bf$ by (see also the notation at the end of Introduction)
$$\F_\Bf^{(j)}:=\F_{\beta^{(1)}} + \cdots + \F_{\beta^{(j)}}, \quad  \F_\Bf:=\F_\Bf^{(m)}, \quad \E_\Bf:=\E/\F_\Bf \quad \text{and} \quad f_\Bf:=f_{\F_\Bf}.$$
We say that $\Bf$ is a {\it leading sequence of powers} for $f$ if 

\medskip\noindent
{\bf (L2)} for $1\leq j\leq m$, the power $\beta^{(j)}$ is leading for the map $f_{\F^{(j-1)}_\Bf}:\Omega\to \E/\F^{(j-1)}_\Bf$;

\medskip\noindent
{\bf (L3)} the smallest closed convex cone in $\R^l$ containing $\beta^{(1)},\ldots,\beta^{(m)}$ and $\R_{\leq 0}^l$, denoted by $\Sigma^-_\Bf$, is salient (i.e., containing no line). In particular, $\Sigma^-_\Bf$ is an unbounded polytope and we have  $\Sigma^-_\Bf\cap \R_{\geq 0}^l=\{0\}$.

\begin{lemma} \label{l:leading-sequence}
Let $\Bf$ be a leading sequence of powers as above. Then the above vector spaces $\F_\Bf^{(j)}$  do not depend on $U$ nor on the chosen coordinate system (S1)-(S2).
\end{lemma}
\proof
We prove the lemma by induction. Assume that  $\F_\Bf^{(j-1)}$  does not depend on $U$ nor on the chosen coordinate system (S1)-(S2).
By applying Lemma \ref{l:dom_power}\'(4) to the map $f_{\F^{(j-1)}_\Bf}:\Omega\to \E/\F^{(j-1)}_\Bf$, we see that $\F_\Bf^{(j)}/\F_\Bf^{(j-1)}$
satisfies the same property. Therefore, the same property holds for  $\F_\Bf^{(j)}$. 
\endproof

Consider a leading sequence of powers $\Bf$ as above. Denote by $\Sigma^{\geq 0}_\Bf$ the smallest closed convex cone in $\R^l$ which contains 
$\R_{\geq 0}^l$ and the family of all powers $\beta$ such that the coefficient of $z'^\beta$ in $f_\Bf$ is non-zero. 
Observe that we obtain the same cone when we only consider minimal $\beta$ in the last family of powers.
Since this set of minimal powers is finite,  the cone $\Sigma^{\geq 0}_\Bf$ is an unbounded polytope. It may be equal to $\R^l$. 

\medskip\noindent
{\bf Complete leading sequences of powers, $\T_\Bf$, $v_{\beta,\Bf}$, $\Bf^0$, $\Bf^+$, $\Sigma_\Bf^0$, $D_\Bf$ and $C_\Bf$.}  We say that a leading sequence of powers $\Bf$ is {\it complete} if it satisfies the following property :

\medskip\noindent
{\bf (L4)} $\Sigma^{\geq 0}_\Bf\cap\Sigma^-_\Bf=\{0\}$.

\medskip

Consider a complete leading sequence $\Bf$ as above. Denote by $\T_\Bf$ the closure of $\pi(\F_\Bf)$ in $\T$ which is a real subtorus of $\T$. 
We will see later that the tori $\T_j$ in Theorem \ref{t:main_high_compact} will be the tori $\T_\Bf$ or their variants.

Denote by $v_{\beta,\Bf}$ the projection of $v_\beta$ in $\E_\Bf$. 
Define $\Sigma^0_\Bf$ as the smallest element (face of smallest dimension) of $\Sigma^{\geq 0}_\Bf$ which contains 0. This is a vector space and we have $\Sigma^0_\Bf=\{0\}$ exactly when $\Sigma^{\geq 0}_\Bf$ is salient. 
Let $\Bf^0$ and $\Bf^+$ be the sets of all $\beta$ in $\Sigma^0_\Bf$ and in $\Sigma^{\geq 0}_\Bf\backslash \Sigma^0_\Bf$, respectively, such that $v_{\beta,\Bf}\not=0$. 
Denote also by $\Sigma_\Bf^+$ the smallest closed convex cone containing $\R_{\geq 0}^l$ and $\Bf^+$. This is also an unbounded polytope. 
The following lemma will be useful for us.

\begin{lemma} \label{l:complete-Hahn-Banach}
Let $\Bf$ be a complete leading sequence of powers as above. Then there is a vector $\lambda=(\lambda_1,\ldots,\lambda_l)\in \Z_{>0}^l$ such that 
$\lambda\cdot\beta^{(j)}<0$ for $1\leq j\leq m$, 
$\lambda\cdot\beta=0$ for $\beta\in \Bf^0$ and $\lambda\cdot\beta >0$ for $\beta\in \Bf^+$. In particular, all powers in $\Bf^0$ are leading for 
$f_\Bf$ and $\Bf^0\setminus\{0\}$ is finite.
\end{lemma}
\proof
Recall that $\Bf, \Bf^0$ and $\Bf^0\cup \Bf^+$ are contained in $\Sigma_\Bf^-, \Sigma_\Bf^0$ and $\Sigma^{\geq 0}_\Bf$ respectively. Moreover, $\Sigma^0_\Bf$ is a vector space and the convex closed cones $\Sigma_\Bf^-, \Sigma_\Bf^{\geq 0}$ are unbounded polytopes. Denote by $\Phi:\R^l\to \R^l/\Sigma^0_\Bf$ the natural projection. 

Observe that the cones $\Phi(\Sigma_\Bf^-)$ and $\Phi(\Sigma_\Bf^{\geq 0})$ are convex, closed, salient cones and also unbounded polytopes.
Then, we can find two disjoint open cones containing  $\Phi(\Sigma_\Bf^-)\backslash\{0\}$ and $\Phi(\Sigma_\Bf^{\geq 0})\backslash\{0\}$.
It follows from Hahn-Banach theorem that there is a linear function $L$ on $\R^l/\Sigma^0_\Bf$ such that $L<0$ on $\Phi(\Sigma_\Bf^-)$ and $L>0$ on $\Phi(\Sigma_\Bf^{\geq 0})\backslash\{0\}$, see \cite[Th.\,I.6]{Brezis}.
By continuity, we can find such a linear form $L$ with rational coefficients.

Consider the linear form $L\circ \Phi$ on $\R^l$. Multiplying it by a suitable positive integer, we can assume that $L\circ \Phi$ is the map $\beta\mapsto\lambda\cdot \beta$ for  some vector $\lambda\in \Z^l$. It is clear that this $\lambda$ satisfies the (in)equalities in the lemma. We also have $\lambda\cdot\beta<0$ for $\beta\in \Sigma_\Bf^-\backslash\{0\}$. Since $\Sigma_\Bf^-$ contains $\R_{\leq 0}^l$, the coefficients of $\lambda$ are strictly positive. 

It remains to prove the last assertion in the lemma. Observe that $\Bf^0\cup\Bf^+$ is exactly the set of powers $\beta$ such that the coefficient $v_{\beta,\Bf}$ of $f_\Bf$ doesn't vanish.  Since the coefficients of $\lambda$ are positive, if $\beta'<\beta$ and $\lambda\cdot\beta=0$ then $\lambda\cdot \beta'<0$. We conclude that all elements of $\Bf^0\setminus\{0\}$ are leading powers for $f_\Bf$. Finally, by Lemma \ref{l:dom_power}\,(1), applied to $f_\Bf$, we obtain that $\Bf^0$ is finite.
\endproof

Consider the vector space $\E_\Bf^{\Bf^0}\simeq \E_\Bf^{|\Bf^0|}$ where $\Bf^0$ is used as a set of indices and $|\Bf^0|$ is the cardinality of $\Bf^0$. More precisely, a vector $u$ in $\E_\Bf^{\Bf^0}$ has $|\Bf^0|$ components  that we denote by $u_\beta$, with $u_\beta\in\E_\Bf$, for $\beta\in\Bf^0$. 
Consider also the natural action of the group $(\C^*)^l$ on $\E_\Bf^{\Bf^0}$ given by 
$$z'\cdot u= \{z'^\beta u_\beta \  : \ \beta\in \Bf^0\} \quad \text{with} \quad z'\in (\C^*)^l \quad \text{and} \quad u=\{u_\beta \  : \ \beta\in \Bf^0\} \in \E_\Bf^{\Bf^0}.$$
The quotient of $\E_\Bf^{\Bf^0}$ by the action of  $(\C^*)^l$ will be denoted by $M$. It is known that $M$ is a complex affine variety and the natural projection $\Phi: \E_\Bf^{\Bf^0} \to M$ is an algebraic map, see Brion \cite[Th.\,1.24]{Brion}. 

We identify $M$ to a Zariski open subset of a projective complex variety $\overline M$. 
Let $w(a)$ denote the image of $\{v_{\beta,\Bf}(a) \ : \ \beta\in \Bf^0\}$ in $M$ for $a\in Z\backslash Y^c_\Rf$. We have the following lemma.

\begin{lemma} \label{l:in-M}
The point $w(a)$ in $M$ does not depend on the coordinate system (S1)-(S2). Moreover, the map $a\mapsto w(a)$ is 
holomorphic from $Z\backslash Y^c_\Rf$ to $M$ and is
meromorphic from  $Z$ to $\overline M$.
\end{lemma}
\proof
The first assertion is a consequence of the formula \eqref{e:coord-change}. The second one is a consequence of the fact that $v_\beta$ is meromorphic on $Z$, holomorphic on $Z\backslash Y^c_\Rf$, and $\Phi$ is algebraic. 
\endproof

Let $\Pi:\E_\Bf^{\Bf^0}\to \E_\Bf$ be defined by 
$$\Pi(u):=\sum_{\beta\in \Bf^0} u_\beta.$$
For every point $w\in M$, define 
$$D_\Bf^w := \big\{\Pi(u) \ : \ u\in\Phi^{-1}(w) \big\}.$$
In particular, we have
$$D_\Bf^{w(a)}=\Big\{\sum_{\beta\in\Bf^0} z'^\beta v_{\beta,\Bf}(a) \ : \ z'\in (\C^*)^l\Big\}.$$ 
The following varieties $D_\Bf$ will involve later in the definition of the algebraic sets $C_j$ in Theorem \ref{t:main_high_compact}.  

\begin{lemma} \label{l:C-B}
Denote by $D_\Bf$ the closure in $\E_\Bf$ of the union of the sets $D_\Bf^{w(a)}$ with $a\in Z\backslash Y^c_\Rf$. Then $D_\Bf$ is an algebraic subvariety of $\E_\Bf$. 
\end{lemma}
\proof
Let  $\M_Z$ denote the set of points $w(a)$ with $a\in Z\backslash Y^c_\Rf.$  This is a subset of $M$ which may not be closed.
We show that $\M_Z$ is algebraic.  
By Lemma \ref{l:in-M}, the map $a\mapsto w(a)$ is meromorphic from $Z$ to $\overline M$. Therefore, $\M_Z$ is a quasi-projective variety in $\overline M$.
Thus, $\Pi(\Phi^{-1}(\M_Z))$ is also a quasi-projective set in $\E_\Bf$. Observe that $D_\Bf$ the closure of the later set in $\E_\Bf$. Therefore, 
$D_\Bf$ should be an algebraic subvariety of $\E_\Bf$. This ends the proof of the lemma.
\endproof

Choose an algebraic subvariety $C_\Bf$ of $\E$ whose projection on $\E_\Bf$ is equal to $D_\Bf$. 
For simplicity, we choose $C_\Bf$ as the intersection of 
$\Pi_{\F_\Bf}^{-1}(D_\Bf)$ with a complement vector space of $\F_\Bf$ in $\E$. Then we have
$$\Pi_{\F_\Bf}^{-1}(D_\Bf)=C_\Bf+\F_\Bf .$$

\smallskip
\noindent
{\bf Construction of complete leading sequences of powers.} We will show now that every cluster value of $\pi(f(z))$ when $z$ tends to $Z\backslash Y^c_\Rf$ can be assigned to a suitable complete leading sequence of powers.

\begin{proposition} \label{p:limit-inclusion}
Let $z_{[1]},z_{[2]},z_{[3]} \ldots$ be a sequence of points in $\Omega\backslash Y$ converging to a point $a\in Z\backslash Y^c_\Rf$. Then, after replacing this sequence by a suitable subsequence, we can find a complete leading sequence of powers $\Bf$ as above such that $f_\Bf(z_{[s]})$ converges to a point in $D^{w(a)}_\Bf$ as $s$ tends to infinity. In particular, when $s$ tends to infinity, the cluster values of $\pi(f_\Bf(z_{[s]}))$ are contained in $\pi(C_\Bf)+\T_\Bf$. 
\end{proposition}

After replacing the sequence $z_{[s]}$ by a suitable subsequence, we can assume that, for every $\beta\in \Z^l$, the sequence $(z_{[j]}')^\beta$ converges to a limit $b_\beta$ which is either 0, $\infty$ or a non-zero complex number. 
Let $\Sigma^{0-}$ denote the smallest convex cone in $\R^l$ containing  all  $\beta\in \Z^l$ with $v_\beta\not=0$ and $b_\beta\not=0,\infty$.  
Denote also by $\Sigma^0$ the minimal element (face of smallest dimension) of $\Sigma^{0-}$ containing $0$.

\begin{lemma} \label{l:leading-Sigma}
We have the following properties.
\begin{enumerate}
\item[(1)] If there is a power $\beta$ with $v_\beta\not=0$ and $b_\beta=\infty$, then there is such a $\beta$ which is a leading power for $f$.
\item[(2)] Assume there is no power $\beta$ with $v_\beta\not=0$ and $b_\beta=\infty$. Then, $\Sigma^{0-}$ is an unbounded  polytope. Moreover, unless $\Sigma^{0-}$ is a vector space, there is a leading power $\beta$ in $\Sigma^{0-}\backslash \Sigma^0$  with $v_\beta\not=0$ and $b_\beta\not=0,\infty$. 
\end{enumerate}
\end{lemma}
\proof
(1) Choose $\beta$ minimal among those with $v_\beta\not=0$ and $b_\beta=\infty$. Clearly, if $\beta'\leq\beta$, then $b_{\beta'}=\infty$.
We see that $\beta$ is minimal among all powers with $v_\beta\not=0$. Moreover, as $b_\beta=\infty$, at least one component of $\beta$ is negative. We conclude that the chosen $\beta$ is a leading power for $f$.

\smallskip

(2) Observe that if $b_\beta\not=0,\infty$, then $\beta$ has a positive component and a negative one, unless $\beta=0$. Moreover, 
 if $b_\beta\not=0,\infty$, then
$b_{\beta'}=0$ for all $\beta'>\beta$ and $b_{\beta'}=\infty$ for all $\beta'<\beta$. 
Consider the set $B$ of all $\beta$ such that $v_\beta\not=0$. We see as in (1) that all $\beta\in B\backslash\{0\}$ with $b_\beta\not=0$ are leading. In particular, 
the set $B^{0-}:=\{\beta\in B \ : \  b_\beta\not=0\}$ is finite and hence 
 $\Sigma^{0-}$ is an unbounded polytope. 

Assume now that $\Sigma^{0-}$ is not a vector space, or equivalently,  $\Sigma^{0-}$ is strictly larger than the vector space $\Sigma^0$. Then the set  $B^{0-}\backslash \Sigma^0$ is non-empty. We have seen that any $\beta$ in $B^{0-}\backslash \Sigma^0$ is a leading power for $f$. This completes the proof of the lemma.
\endproof

We continue the proof of Proposition \ref{p:limit-inclusion}. 
We need to construct a sequence $\Bf$ as above satisfying (L2)-(L4) and the properties stated in Proposition \ref{p:limit-inclusion}. This sequence $\Bf$ may be empty.
We will use the notations introduced above, in particular, the ones introduced near the conditions (L1)-(L4).
The construction is done by induction as follows. 

If there is a $\beta$ such that $v_\beta\not=0$ and $b_\beta=\infty$, by Lemma \ref{l:leading-Sigma}\,(1), we choose $\beta^{(1)}$ a leading power for $f$ with $v_{\beta^{(1)}}\not=0$ and $b_{\beta^{(1)}}=\infty$. Otherwise, we apply Lemma \ref{l:leading-Sigma}\,(2) and choose $\beta^{(1)}$ a leading power for $f$ with $v_{\beta^{(1)}}\not=0$ and $b_{\beta^{(1)}}\not=0,\infty$ that belongs to $\Sigma^{0-}\backslash \Sigma^0$, unless the later set is empty. When  $\Sigma^{0-}\backslash \Sigma^0$ is empty,  $\Bf$ is chosen to be empty too.

Assume that $\beta^{(1)}, \ldots,\beta^{(s)}$ are constructed. We construct $\beta^{(s+1)}$ in the same way as we did for $\beta^{(1)}$. More precisely, we replace $f$ by $f_{\F^{(s)}_\Bf}$ and apply Lemma \ref{l:leading-Sigma} to this map instead of $f$. We end the construction when we reach the case where $\Sigma^{0-}\backslash \Sigma^0$ is empty. In particular, we end the construction if we have $\Sigma^{0-}=\{0\}$. 
The construction gives us a  sequence $\Bf=\{\beta^{(1)},\ldots, \beta^{(m)}\}$ satisfying (L2). 

In order to simplify the notation, we define for every $\beta\in \R^l$
$$|b_\beta|:=\lim_{s\to\infty} |z'_{[s]}|^\beta$$
if this limit (finite or infinite) exists.

\begin{lemma} \label{l:b_beta}
The following properties hold.
\begin{itemize}
\item[(1)]  $|b_\beta|$ exists and $|b_\beta|\not=0$ for every $\beta\in \Sigma^-_\Bf$.
\item[(2)]  $|b_\beta|$ exists and  $|b_\beta|\not=0,\infty$ for every $\beta\in \Sigma^0_\Bf$.
\item[(3)]  $|b_\beta|$ exists and  $|b_\beta|=0$ for every $\beta\in \Sigma^{\geq 0}_\Bf\backslash \Sigma^0_\Bf$.
\end{itemize}
\end{lemma}
\proof
(1) Consider $\beta\in \Sigma^-_\Bf$. By definition of $\Sigma^-_\Bf$, there is $\beta'\leq 0$ and non-negative numbers $\lambda_j$ such that 
$$\beta=\beta'+\sum_j \lambda_j \beta^{(j)} \quad \text{and hence} \quad |z'|^\beta = |z'|^{\beta'}\prod_j |z'|^{\lambda_j \beta^{(j)}}.$$
Since $b_{\beta^{(j)}}\not=0$, we easily deduce that $|b_\beta|$ exists and not equal to 0.

\smallskip

(2) Observe that there is a finite set $B^0$ of  
$\beta'\in\Sigma^0_\Bf\cap \Z^l$ with $v_{\beta',\Bf}\not=0$ such that $B^0$ spans $\Sigma^0_\Bf$ and the convex hull of $B^0$ contains 0 in its interior. 
So we can find positive numbers $\lambda_{\beta'}$ such that
$$\sum_{\beta'\in B^0} \lambda_{\beta'} \beta'=0 \quad \text{hence} \quad \prod_{\beta'\in B^0} |z'|^{\lambda_{\beta'} \beta'}=1 \quad \text{and} \quad \prod_{\beta'\in B^0} |b_{\beta'}|^{\lambda_{\beta'}}=1.$$
In particular, we have $|b_{\beta'}|\not=0,\infty$ for $\beta'\in B^0$. Since any $\beta\in \Sigma^0_\Bf$ can be write as a linear combination of $\beta'\in B^0$, we deduce that $|b_\beta|$ exists and is not equal to 0 nor $\infty$. 

\smallskip

(3) As in (1), we obtain that $|b_\beta|\not=\infty$ for $\beta\in \Sigma^{\geq 0}_\Bf$. 
By the construction of $\Bf,$ we have $b_\beta=0$ for every $\beta\in \Z^l$ with $v_{\beta,\Bf}\not=0$ unless $\beta\in \Sigma^0_\Bf$. Using the same argument as in (1) or (2), we easily deduce that 
$|b_\beta|=0$ for every $\beta\in \Sigma^{\geq 0}_\Bf\backslash \Sigma^0_\Bf$.
\endproof

\proof[{\bf End of the proof of Proposition \ref{p:limit-inclusion}}]
We show that the set of powers $\Bf$ satisfies (L3).  Assume by contradiction that this is not true. Then there are non-negative numbers $\lambda_j$ and a vector $\beta\in\R_{\leq 0}^l$, not all equal to 0, such that 
$$\beta+\sum_j \lambda_j \beta^{(j)} =0 \quad \text{and hence} \quad |b_\beta|\prod_j |b_{\beta^{(j)}}|^{\lambda_j}=1.$$
This and Lemma \ref{l:b_beta}\,(1) imply $|b_\beta| \not = \infty$ and $ |b_{\beta^{(j)}}|\not=\infty$ when $\lambda_j\not=0$. 
We then deduce that $\beta=0$. If $J$ is the set of indices $j$ such that $\lambda_j\not=0$, we obtain that $|b_{\beta^{(j)}}|\not=\infty$ for $j\in J$ and the convex hull of $\beta^{(j)}$ with $j\in J$ contains 0. From the construction of $\Bf$, we see that $\beta^{(j)}\in \Sigma^{0-}\backslash\Sigma^0$ for $j\in J$ and their convex hull cannot contain 0. This is a contradiction. So 
the set $\Bf$ satisfies (L3).         

We show that $\Bf$ satisfies (L4). By Lemma \ref{l:b_beta}\,(1)(3), the intersection $\Sigma_\Bf^-\cap\Sigma_\Bf^{\geq 0}$ is contained in $\Sigma^0_\Bf$. Observe that the later set is a vector space contained in the vector space $\Sigma^0$ in the construction of $\Bf$. As in the proof of (L3), we obtain that $\Sigma_\Bf^-\cap\Sigma^0=\{0\}$. So we also have  $\Sigma_\Bf^-\cap\Sigma_\Bf^{\geq 0}=\{0\}$ which is exactly the property (L4). 

Finally, by Lemma \ref{l:b_beta}\,(2), we have $b_\beta\not=0$ for $\beta\in\Bf^0$. It is not difficult to see that 
$$v:=\sum_{\beta\in \Bf^0} b_\beta v_{\beta,\Bf}(a)$$
belongs to $D^{w(a)}_\Bf$. We conclude that $f_\Bf(z_{[s]})$ converges to the point $v$ in $D^{w(a)}_\Bf$ as $s$ tends to infinity.

For the last assertion in the proposition, it is enough to observe that the distance between $f_\Bf(z_{[s]})$ and $v+\F_\Bf$ tends to 0. The result follows easily.
\endproof

\noindent
{\bf Construction of a good holomorphic disc.} Using Proposition \ref{p:limit-inclusion}, in order to obtain Theorem \ref{t:main_high_compact}, we still need to show that all points in $\pi(C_\Bf)+\T_\Bf$ are cluster values of $\pi(f_\Bf(z))$ when $z$ tends to $Z\backslash Y^c_\Bf$. For this purpose, 
as already mentioned in Introduction, we will construct a good holomorphic disc in $\Omega$ in order to apply Theorem \ref{t:compact-1D}.

\begin{proposition} \label{p:good-disc}
Let $\Bf$ be a complete leading sequence of powers as above. Then, for every point $a\in Z\backslash Y^c_\Rf$ and $b\in D_\Bf^{w(a)}$, there is a holomorphic map $\phi:\D\to \Omega$ such that
\begin{enumerate}
\item[(1)] $\phi(0)=a$ and $\phi(\D^*)\subset \Omega\backslash Y$;
\item[(2)] $\F_\Bf$ is the smallest vector subspace of $\E$ such that $(f\circ\phi)_{\F_\Bf}$ extends holomorphically through $0\in\D$;
\item[(3)] $(f\circ\phi)_{\F_\Bf}(0)=b$.
\end{enumerate} 
In particular, the set $\pi(C_\Bf)+\T_\Bf$ is contained in $\Lambda_{f,Z}$. 
\end{proposition}

Observe that if there is a holomorphic disc satisfying (1)(2)(3) as above, by Theorem \ref{t:compact-1D}, the set of cluster values of $\pi(f(\phi(x))$, when $x$ tends to 0, is equal to $\pi(b)+\T_\Bf$. So, the last assertion in the proposition is a consequence of the previous one. In order to complete the proof of the proposition, we need the following lemmas.

\begin{lemma} \label{l:independence}
Let $h_1,\ldots, h_p$ be holomorphic functions in a neighbourhood of $0\in\C^l$ which are linearly independent. Denote by $h_{1,N},\ldots, h_{p,N}$ the sums of terms of degree at most $N$ in their Taylor expansions. Then for every $N$ large enough  the functions $h_{1,N},\ldots,h_{p,N}$ are linearly independent.
\end{lemma}
\proof
Write 
$$h_j(z')=\sum_{\beta\in\Z_{\geq 0}^l} a_{j,\beta} z'^{\beta} \quad \text{with} \quad a_{j,\beta}\in\C$$
the Taylor expansion of $h_j$. Let $r$ denote the rank of the matrix 
$$(a_{j,\beta})_{1\leq j\leq p, \beta\in\Z_{\geq 0}^l}.$$
Fix an integer $N_0$ large enough such that for $N\geq N_0$  the rank of the matrix
$$A_N:=(a_{j,\beta})_{1\leq j\leq p, |\beta|\leq N}$$
is maximal, i.e., equal to $r$. 

The lemma is equivalent to the equality $r=p$. 
Assume by contradiction that $r<p$. Observe that the space $V_N$ of all vectors $v\in\C^p$ such that $A_Nv=0$ is of dimension $p-r$ for all $N\geq N_0$. Since $V_N$ decreases when $N$ increases, we deduce that $V_N=V_{N_0}$ for $N\geq N_0$. If $v=(v_1,\ldots, v_p)$ is a non-zero vector in $V$, we see that $A_\infty v=0$ and hence $v_1h_1+\cdots+v_ph_p=0$. This is a contradiction and ends the proof of the lemma.
\endproof

\begin{lemma} \label{l:power-separation}
For every integer $N\geq 1$, there are a holomorphic monomial map
$\psi:\D\to \C^l$ with $\psi(x)=(x^{\gamma_1},\ldots,x^{\gamma_l})$ and $\gamma_j\in\Z_{>0}$, and an integer $M\geq 1$ such that 
\begin{enumerate}
\item[(1)]  the monomials $\psi(x)^\beta$ are distinct for $\beta\in\Z_{\geq 0}^l$ and $|\beta|\leq N$; 
\item[(2)] the degree of $\psi(x)^\beta$ is less than $M$ for $\beta\in\Z_{\geq 0}^l$ with $|\beta|\leq N$ and at least $M$ otherwise.
\end{enumerate}
\end{lemma}
\proof
Choose a vector $\gamma=(\gamma_1,\ldots,\gamma_l)$ in $\R_{>0}^l$ such that $\gamma^*\leq \gamma_j < (1+N^{-1})\gamma^*$ for every $1\leq j\leq l$ and some $\gamma^*>0$, and  the map $\beta\mapsto \gamma\cdot\beta$ is injective on the set $\{\beta\in \Z_{\geq 0}^l \ : \ |\beta|\leq   N \}$. By continuity, we can choose such a $\gamma$ with rational coefficients. Multiplying this vector by a suitable positive integer allows us to assume that $\gamma\in\Z_{>0}^l$. We also replace $\gamma^*$ by the minimum of the $\gamma_j$'s which is an integer number.
Define $M:=(N+1)\gamma^*$ and $\psi(x) :=(x^{\gamma_1},\ldots,x^{\gamma_l})$. It is clear that this choice satisfies the lemma.
\endproof

\proof[{\bf Proof of Proposition \ref{p:good-disc}}]
Fix a coordinate system $z=(z',z'')$ as above which is centered at the point $a$. Each coefficient $v_\beta$ will be considered as a holomorphic function in $z''$ with values in $\E$.
We can choose a point $\alpha\in (\C^*)^l$ such that 
$$\sum_{\beta\in\Bf^0} \alpha^\beta v_{\beta,\Bf}(0)=b.$$
Let $\lambda$ be the vector in Lemma \ref{l:complete-Hahn-Banach}.  Let $\Bf_\lambda$ denote the set of powers $\beta$ such that $\lambda\cdot\beta<0$ and $v_\beta\not=0$.  This is a finite set and the space $\F_\Bf$ is contained in the vector space spanned by $v_\beta(z'')$ with $\beta\in\Bf_\lambda$ and $z''$ in a neighbourhood of $0\in \C^{q-l}$. Actually, these vector spaces are equal but we don't need this stronger property.

Consider the Taylor expansion
\begin{equation} \label{e:v_beta}
v_\beta(z'') =\sum_{\vartheta\in\Z_{\geq 0}^{q-l}} z''^\vartheta v_{\beta,\vartheta} \quad \text{with}\quad  v_{\beta,\vartheta}\in\E.
\end{equation}
Observe that if  $N_0$ is a large enough integer, then for every $\beta\in\Bf_\lambda$ and every $z''$ in a neighbourhood of $0\in \C^{q-l}$, the vector $v_\beta(z'')$ belongs to the space spanned by the vectors $v_{\beta,\vartheta}$ with  $|\vartheta|\leq N_0$. We use here the fact that $\Bf_\lambda$ is finite.

Observe also that the functions $z'^\beta$ with $\beta\in\Bf_\lambda$ are linearly independent. Therefore, the functions $(\alpha+z')^\beta$ are also linearly independent. Consider the Taylor expansion 
$$(\alpha+z')^\beta=\sum_{\theta\in\Z_{\geq 0}^l} c_{\beta,\theta} z'^\theta \quad \text{with}\quad c_{\beta,\theta}\in\C.$$
By Lemma \ref{l:independence}, we can fix an integer $N_0$ large enough so that the polynomials 
$$\sum_{|\theta|\leq N_0} c_{\beta,\theta} z'^\theta \quad \text{with} \quad \beta\in\Bf_\lambda$$ 
are linearly independent, or equivalently, the matrix 
$$(c_{\beta,\theta})_{\beta\in\Bf_\lambda,|\theta|\leq N_0}$$
is of maximal rank $|\Bf_\lambda|$ (the cardinality of $\Bf_\lambda$).

Consider the vector space $\C^{\Bf_\lambda}\simeq \C^{|\Bf_\lambda|}$ where we use $\Bf_\lambda$ as a set of indices. 
We deduce from the last discussion that the family of vectors
\begin{equation} \label{e:vectors}
(c_{\beta,\theta})_{\beta\in\Bf_\lambda} \quad \text{with} \quad |\theta|\leq N_0
\end{equation}
span the space $\C^{\Bf_\lambda}$.

Define $N:=2N_0$. Let $\psi$ and $M$ be as in Lemma \ref{l:power-separation}.  Define the map $\phi:\D\to \Omega$ by 
$$\phi(x):=\big(x^{M\lambda_1}(\alpha_1+\psi_1(x)),\ldots,x^{M\lambda_l}(\alpha_l+\psi_l(x)), \psi_{l+1}(x),\ldots,\psi_q(x)  \big).$$
Consider a vector subspace $\F$ of $\E$ such that $(f\circ\phi)_\F$ is bounded. We only need to show that $\F$ contains $\F_\Bf$.  For this goal, we are interested in the non-zero terms with negative power in the Laurent expansion of the map $f_\F(\phi(x))$. 

Assume by contradiction that $\F$ doesn't contain $\F_\Bf$.
Then, the family $\Bf_\lambda'$ of all powers $\beta\in\Bf_\lambda$ such that  
we don't have $v_\beta(z'')$ contained in $\F$ for all $z''$ near 0, is non-empty.
Consider the family $\Bf_\lambda^{\min}$ of all powers $\beta\in\Bf'_\lambda$ such that  
$\lambda\cdot \beta$ is minimal. Denote by $-A$ this minimal value of $\lambda\cdot \beta$ which is a negative integer. 

We use now Properties (1), (2) in Lemma \ref{l:power-separation} and the expansions \eqref{e:Laurent-bis}, \eqref{e:v_beta}. We have that 
the partial sum of terms of degree less than $-AM+M$ in $f_\F(\phi(x))$ is equal to
\begin{equation} \label{e:partial-sum}
\sum_{|\theta|+|\vartheta|\leq N \atop \beta\in \Bf_\lambda^{\min}}  x^{-AM} \psi(x)^{(\theta,\vartheta)}  c_{\beta,\theta} v_{\beta,\vartheta}.
\end{equation}
Here, $(\theta,\beta)$ is an element of $\Z_{\geq 0}^q$, $\psi(x)^{(\theta,\vartheta)}$ is a power of $x$ of degree less than $M$, and for simplicity, we still denote by $v_{\beta,\vartheta}$ its image in $\E/\F$. 
Indeed, if we consider $\beta\not\in \Bf_\lambda^{\min}$, then the involving powers of $x$  in $f_\F(\phi(x))$  are at least equal to $(-A+1)M$, and if we consider $(\theta,\beta)$  with $|\theta|+|\vartheta|\geq N+1$ then $\psi(x)^{(\theta,\vartheta)}$ is a power of $x$ of degree at least $M$.

Observe now that the terms in  \eqref{e:partial-sum} contain different negative powers of $x$. 
According to Theorem \ref{t:compact-1D}, the space $\F$ contains the vector 
$$\sum_{\beta\in \Bf_\lambda^{\min}} c_{\beta,\theta}  v_{\beta,\vartheta} $$
for all $\theta, \vartheta$ such that $|\theta|+|\vartheta|\leq N$, in particular, for $|\theta|\leq N_0$ and $|\vartheta|\leq N_0$. Since the vectors in \eqref{e:vectors} span $\C^{\Bf_\lambda}$,  
we deduce that $v_{\beta,\vartheta}$ belongs to $\F$ for all $\beta\in \Bf_\lambda^{\min}$ and $|\vartheta|\leq N_0$. 
Hence, by the choice of $N_0$, the vector $v_\beta(z'')$ belongs to $\F$ for all $\beta\in \Bf_\lambda^{\min}$ and $z''$ small enough.
This contradicts the definition of $\Bf_\lambda^{\min}$ and ends the proof of the proposition.
\endproof

\proof[{\bf Proof of Theorem \ref{t:main_high_compact}}]
Consider all manifolds $Z$ which are contained in $Y^*$ and all complete leading sequences of powers $\Bf$ as above, together with the associated vectors spaces, tori and algebraic varieties $\F_\Bf$, $\T_\Bf$, $D_\Bf$, $C_\Bf$. By Propositions \ref{p:limit-inclusion} and \ref{p:good-disc}, we have 
\begin{equation} \label{e:union}
\Lambda_{f,Y^*}=\bigcup_\Bf  (\pi(C_\Bf) + \T_\Bf).
\end{equation}
So we obtain Theorem \ref{t:main_high_compact} except the last assertion. In particular, the proof of Theorem \ref{t:main-2} is now complete.

We will call $(\F_\Bf,D_\Bf)$ {\it a basic pair} of our construction.
In order to get the last assertion in Theorem \ref{t:main_high_compact}, we will add to the above family of $(\F_\Bf,D_\Bf)$ other pairs $(\F_j,D_j)$ without changing the union on the right-hand side of \eqref{e:union}. 

Consider a basic pair $(\F_\Bf,D_\Bf)$ such that $D_\Bf$ is of positive dimension. We can apply our construction to the inclusion map from $D_\Bf$ to $\E_\Bf$ and get basic pairs for this map. Let $(\F'_j,D_j)$ denote one of them. Define $\F_j:=(\Pi_{\F_\Bf})^{-1}(\F_j')$, $\T_j:=\overline{\pi(\F_j)}$ and $C_j$ the intersection between $\Pi_{\F_j}^{-1}(D_j)$ with a complement vector space of $\F_j$ in $\E$. It is not difficult to see that $\pi(C_j)+\T_j$ is contained in the closure of $\pi(C_\Bf) + \T_\Bf$. As $\Lambda_{f,Y^*}$ is closed, $\pi(C_j)+\T_j$ is contained in $\Lambda_{f,Y^*}$. 

We add these pairs $(\F_j,D_j)$ to the family of basic pairs constructed earlier. This doesn't change the union  on the right-hand side of \eqref{e:union}.
The new pairs will be called {\it pairs of second generation}. It is clear that if $D_\Bf$ is of positive dimension, then $\F_\Bf$ is strictly contained in the space $\F_j$ from at least one pair of second generation. 

In the same way, we repeat the construction for the new pairs and get pairs of third generation. We repeat the construction until we get a generation  with some $\F_j$ equal to $\E$, or otherwise,  all $D_j$ finite. In the first case, we always have $\Lambda_{f,Y^*}=\T$ and we can just replace $C_j$ by any point of $\E$ in order to get the result. In the second case, it is clear that when $\F_j$ is maximal for the inclusion, then $D_j$ is finite and hence $C_j$ is finite. This ends the proof of Theorem \ref{t:main_high_compact}.
\endproof

\noindent
{\bf Example of algebraic flow on a semi-torus.} We will end this section with a concrete example showing that the case of a non-compact semi-torus is of another nature.
Define $\E:=\C^3$ and $\T:= (\C/\Z)^3$.
Consider the algebraic hypersurface $X$ of  $\C^3$ which is the image of the holomorphic map $\tau: \C^2 \to \C^3$  given by 
$$\tau(z):= (e^{i\pi/4}z_1, z_2, z_1 z^2_2+ z_1 z_2) \quad \text{with} \quad z=(z_1,z_2)\in\C^2.$$
Clearly, $X$ is closed in $\C^3$ as $\tau(z)$ tends to infinity as $z$ tends to infinity.

Denote by $\Lambda_0$ the set of limit points of $\pi(\tau(z))$ as  $z_1 \to \infty$ and $z_2 \to \infty.$
Denote also by $\Lambda_1$ (resp. $\Lambda_2$) the set of limit points of $\pi(\tau(z))$ when  $z_1 \to \infty$ and $z_2$ is  bounded (resp. $z_2$ tends to infinity and $z_1$ is bounded).   Hence,  $\Lambda_0 \cup \Lambda_1\cup \Lambda_2$ is the set of all limit points of $\pi(\tau(z))$ when $z$ tends to infinity.

Define $\sigma(z_2):= (0,z_2,0)$ and consider the set 
$$C_1:=\sigma(L) \quad \text{with} \quad  L:=\big\{s \in \C: \arg(s^2+s)=\pi/4 \text{ or } -3\pi/4\big\} \cup \big\{s\in \C:\  s^2+s=0\big\}.$$
Observe that $L$ is the pullbacks of the real line $\arg(s)=\pi/4, -3\pi/4$ (we include here the point $s=0$) by the map $s\mapsto s^2+s$. This map has a unique critical value $-1/4$ which is outside the considered real line.
So $L$ is a union of two disjoint real curves which are closed subsets of $\C$.
So $C_1$ is also a union of two disjoint real curves which are closed subsets of $\{0\}\times\C\times\{0\}$ in $\C^3$.

Finally, define
$$\F_1:=\C\times\{0\}\times\C \quad \text{and} \quad  \quad  \F_2:=(e^{i \pi/4}\R)\times \C^2.$$
Then, consider 
$$\T_1:= \F_1/ (\F_1\cap \Z^3) \quad \text{and} \quad \T_2:= \F_2/(\F_2 \cap \Z^3).$$
Both of them are closed in $\T$, $\dim_\C \T_1=\dim_\C\F_1=2$ and $\dim_\R \T_2=\dim_\R \F_2=5$.

\medskip\noindent
{\bf Claim. } We have
$$\Lambda_0= \varnothing, \quad  \Lambda_1= \pi(C_1)+ \T_1  \quad \text{and} \quad 
\Lambda_2= \T_2 . $$

\proof[{\bf Proof of the claim}]
In order to describe the cluster points of $\pi(f(z))$ when $z$ tends to infinity, we only need to consider the situation where $z$ tends to infinity but $\pi(\tau(z))$ stays in a compact subset of $\T$ or equivalently
\begin{equation*} \tag{$\dagger$}
|\Im(e^{i\pi/4} z_1)|, \ |\Im(z_2)|,\ |\Im(z_1z_2^2+z_1z_2)| \ \ \text{are bounded by a constant}. 
\end{equation*}

\smallskip\noindent
{\bf Case 1.} 
Assume that  $z_1 \to \infty$ and $z_2 \to \infty.$ We deduce from ($\dagger$) that 
both  $\arg(e^{i\pi/4} z_1)$ and $\arg(z_2)$ tend to $0$ or $\pi$. So
 $\arg z_1 \to -\pi/4$ or $3\pi/4$ and $\arg z_2 \to 0$ or $\pi.$  By writing $z_1 z_2^2 +z_1 z_2=z_1z_2^2(1+1/z_2)$, we see that $\Im(z_1 z_2^2 +z_1 z_2)$ can not be bounded. Consequently, we get $\Lambda_0=\varnothing.$

\medskip\noindent
{\bf Case 2.} Assume that $z_1$ tends to infinity and $z_2$ is bounded. So we can assume that $z_2$ converges to a point $s\in\C$.
We deduce from ($\dagger$) that $\arg(e^{i\pi/4} z_1)$ tends to $0$ or $\pi$ and $\arg(z_1z_2^2+z_1z_2)$ tends to 0 or $\pi$ unless $s^2+s=0$.
Since $\arg(z_1z_2^2+z_1z_2)=\arg(z_1)+\arg(z_2^2+z_2)$,  we deduce that $s$ belongs to $C_1$.
Hence, $\Lambda_1$ is contained in $\pi(C_1)+ \T_1$. 

In order to obtain the second identity in the claim, consider three arbitrary points $s\in C_1$, $a_1\in \C$ and $a_3\in\C$ . 
We need to check that $\pi(a_1,s,a_3)$ belongs to $\Lambda_1$. 
We only consider the case where $\arg(s^2+s)=\pi/4$ because the case with $\arg(s^2+s)=-3\pi/4$ or $s^2+s=0$ can be treated in the same way. 

Consider a large integer parameter $N>0$ and choose $(z_1,z_2)$ as  a solution of the following equations
$$e^{i\pi/4} z_1 = a_1+ N  \quad \text{and} \quad  z_1z_2^2+z_1z_2 =a_3+ \left \lceil e^{-i\pi/4}N(s^2+s) \right \rceil .$$
There is a unique choice for $z_1$ and we have $z_1=e^{-i\pi/4}N+O(1)$ as $N$ tends to infinity. There are 
two choices for $z_2$ but we choose $z_2$ so that $z_2=s+o(1)$ as $N$ tends to infinity. Clearly, when $N$ goes to infinity, $z_1$ tends to infinity, $z_2$ tends to $s$ and $\pi(\tau(z))$ tends to $\pi(a_1,s,a_3).$ This completes the proof of the second identity in the claim.

\medskip\noindent
{\bf Case 3.} Assume that $z_2$ tends to infinity and  $z_1$ is bounded.
So we can assume that $z_1$ converges to a point $a_1\in\C$. If $a_1=0$, then clearly all obtained cluster values of $\pi(\tau(z))$ are in $\T_2$.
Otherwise, we have $a_1\not=0$. Then as above,  $\arg(z_2)$ tends to 0 or $\pi$ and 
$\arg(z_1z_2^2+z_1z_2)\simeq \arg(a_1)+2\arg(z_2)$ also tends to $0$ or $\pi$. Hence $\arg a_1=0$ or $\pi$. Thus, all obtained cluster values of $\pi(\tau(z))$ belong to $\T_2$ as well. We have shown that $\Lambda_2\subset \T_2$.

It remains to check that $\T_2\subset\Lambda_2$. Consider any point $(t,b_2,b_3)\in \R\times\C^2$. We need to show that $\pi(e^{i\pi/4}t,b_2,b_3)$ belongs to $\Lambda_2$. Consider a large integer parameter $N>0$ and $(z_1,z_2)$ the unique solution of the equations
$$z_2=b_2+N \quad \text{and} \quad z_1(z_2^2+z_2) = b_3 + \left \lceil tN^2 \right \rceil .$$
Clearly, $z_1$ tends to $t$ and $\pi(\tau(z))$ tends to $\pi(e^{i\pi/4}t,b_2,b_3)$ as $N$ goes to infinity. This ends the proof of the claim.
\endproof

\section{Ax-Lindemann-Weierstrass theorem} \label{s:ALW}

In this section, we will prove Theorem \ref{t:ALW}. Let $\E, \T:=\E/\Gamma$ and $\pi$ be as above, where $\T$ is not necessarily compact. 
We need the following basic lemma.

\begin{lemma} \label{l:inv-Lie} Let $\H$ be a connected real subgroup of $\T$ and let $\overline \H^\Zariski$ be the Zariski closure of $\H$ in $\T$. Then 
$\overline \H^\Zariski$ is the smallest closed complex Lie subgroup of $\T$ which contains $\H$.  Moreover, if 
$W$ is a real vector subspace of $\E$ which is contained in the complex vector space spanned by  $W\cap \Gamma_\R$ and  
$\H=\pi(W/(W \cap \Gamma))$, then $\overline \H^\Zariski$ is a sub-semi-torus of $\T.$
\end{lemma}
\proof
For the first assertion, it is enough to check that $\G:=\overline \H^\Zariski$ is a group.
Consider the map $\Pi(a,b):= a-b$ from $\T^2$ to $\T$. We only have to show that the image of $\G^2$ by $\Pi$ is contained in $\G$.
Observe that the Zariski closure of $\H^2$ in $\T^2$ contains all sets $\{a\}\times\G$ with $a\in \H$. So it contains $\H\times\{b\}$ for all $b\in\G$. Hence, it contains the Zariski closure of the latter set which is $\G\times\{b\}$. We conclude that the Zariski closure of $\H^2$ in $\T^2$ is equal to $\G^2$.

On the other hand, since $\H$ is a subgroup of $\T,$ the set $\Pi(\H^2)$ is contained in $\H.$ This coupled with the fact that $\H \subset \G$ shows that $\Pi^{-1}(\G)$ is a Zariski closed subset of $\T^2$ containing $\H^2.$ It follows that $\Pi^{-1}(\G)$ contains $\G^2$, or equivalently, $\Pi(\G^2)$ is contained in $\G$.  Hence the first assertion in the lemma follows.

We prove now the second  assertion. 
Let $V$ be the complex vector subspace of $\E$ such that $\pi(V)=\overline\H^\Zariski$. Observe that $\pi(V\cap\Gamma_\R)=\overline\H^\Zariski\cap \pi(\Gamma_\R)$ is compact because $\overline \H^\Zariski$ is closed in $\T$ and $\pi(\Gamma_\R)=\Gamma_\R/\Gamma$ is compact. Let $W_1$ be the complex vector space spanned by $V\cap\Gamma_\R.$ Since $V$ is complex, we have 
$W_1\subset V$. Therefore, $W_1\cap\Gamma_\R= V\cap\Gamma_\R$, 
$\pi(W_1\cap\Gamma_\R)=\overline\H^\Zariski\cap \pi(\Gamma_\R)$ and $W \subset W_1$ by hypothesis on $W$. In particular, $\pi(W_1\cap\Gamma_\R)$ is compact and $\H':=\pi(W_1)$ is a semi-torus contained in $\overline\H^\Zariski=\pi(V)$. We only need to check that $\H'$ is closed in $\T$ because this property together with the inclusion $W \subset W_1$ implies that $\overline \H^\Zariski=\H'$. 

Observe that $\H'\cap \pi(\Gamma_\R)$ is compact because it is equal to $\pi(W_1\cap\Gamma_\R)$.  Consider the natural projection $P:\T\to \T/\pi(\Gamma_\R)$ which is a proper map. Note that $\T/\pi(\Gamma_\R) = \E/ \Gamma_\R.$ So the composition map $P \circ \pi: \E \to \T/ \pi(\Gamma_\R)$ is exactly the projection from $\E$ to $\E/ \Gamma_\R.$ Since $P(\H')= P \circ \pi (W_1),$  $P(\H')$ is simply a vector subspace of  $\E/ \Gamma_\R$ which is obviously closed in $\E/ \Gamma_\R.$ Hence $\H'$ is closed in $\T$. This ends the proof of the lemma.
\endproof

The following result is the first main theorem in this section.

\begin{theorem} \label{t:ALW-1D}  
Under the hypotheses (H0)-(H3) as in Theorem \ref{t:1D-Gamma}, denote by $\F$ the smallest complex subspace of $\E$ such that the map $f_\F:\D^*\to\E/\F$ extends holomorphically through $0$, see also Lemma \ref{l:F-l}. Then 
the Zariski closure $\overline{\pi(f(\D^*))}^\Zariski$ of   $\pi(f(\D^*))$ in $\T$ is invariant under the action of $\T_\F:=\overline{\pi(\F)}^\Zariski$. In particular, if $W$ denotes the Zariski closure of the image of $f(\D)$ in $\T/\T_\F$, then $\overline{\pi(f(\D^*))}^\Zariski$ is the pullback of $W$ by the natural projection from $\T$ to $\T/\T_\F$. 
\end{theorem}
\proof
Note that since $\T_F$ is a closed subgroup of $\T$, the quotient $\T/\T_\F$ is a commutative complex Lie group. 
Denote by $p:\T\to\T/\T_\F$ the natural projection.  So the map $p\circ\pi\circ f$ extends to a holomorphic map from $\D$ to $\T/\T_\F$.
Without loss of generality, we can assume for simplicity that $p(\pi(f(0)))=0$. Define $S:= p(\pi(f(\D)))$. So either $S$ is a Riemann surface or it is the singleton $\{0\}$.

Consider the case where $S$ is a Riemann surface. By Lemma \ref{l:F-l-2} and Theorem \ref{t:1D-Gamma}, there is a real vector subspace $\F'$ of $\E$ such that $\F'+i\F'=\F$ and the usual closure of $\pi(f(\D^*))$ in $\T$ contains $\pi(\F')$. It follows that the Zariski closure of $\pi(f(\D^*))$ in $p^{-1}(S)$, denoted by $K$, contains $\pi(\F')$. The fact that $\F'+i\F'=\F$ implies that 
$\pi(\F')$ is Zariski dense in $\pi(\F)$. So $K$ contains $\pi(\F)$ and hence $\T_\F$ because $\pi(\F)$ is Zariski dense in $\T_F$. Since $\pi(f(\D^*))$ is irreducible, $K$ is also irreducible. Then, in the variety $p^{-1}(S)$, 
the fact that $K$ contains the hypersurface $\T_\F$ together with $\pi(f(\D^*))$ implies that $K=p^{-1}(S)$. 

It is not difficult to see that the last identity still holds when $S=\{0\}$. 
We deduce that, in any case, $\overline{\pi(f(\D^*))}^\Zariski$, which contains $K$, contains $p^{-1}(S)$ as well. So $\overline{\pi(f(\D^*))}^\Zariski$ is the Zariski closure of $p^{-1}(S)$ in $\T$.
Finally, if $a$ is a point in $\T_\F$ then $a+\overline{\pi(f(\D^*))}^\Zariski$ is Zariski closed in $\T$ and contains $p^{-1}(S)$ as the latter set is invariant by $\T_F$. It follows that $a+\overline{\pi(f(\D^*))}^\Zariski$ is equal to $\overline{\pi(f(\D^*))}^\Zariski$. Thus, we obtain the first assertion in the theorem. The second assertion follows easily.
\endproof

\proof[{\bf Proof of Theorem \ref{t:ALW}}] 
Let $\T'$ be the maximal sub-semi-torus of $\T$ such that $\overline{\pi(X)}^\Zariski$ is invariant by $\T'$. Let $\F'$ be the complex subspace of $\E$ such that $\pi(\F')=\T'$. Notice that since the projection of $X$ on $\E/\F'$ is an algebraic variety, it is a dense Zariski open subset of its closure.

Replacing $\E,\T$ and $X$ by $\E/\F', \T/\T'$ and the closure of the projection of $X$ on $\E/\F'$, we can assume that $\overline{\pi(X)}^\Zariski$ is not invariant by any non-trivial sub-semi-torus of $\T$. 
In order to get Theorem \ref{t:ALW}, it is enough to prove the following claim.

\medskip\noindent
{\bf Claim.} The restriction of $\pi$ to $X$ is a proper map from $X$ to $\T$. In particular, if $\pi(X)$ is relatively compact in $\T$, then $X$ is a point.

\medskip

Note that if $\pi(X)$ is relatively compact in $\T$, then the first assertion in the claim implies that $X$ is compact and therefore, $X$ is a point. So we only need to prove the first assertion in the claim.

Consider $\E$ as an affine chart of the projective space $\P(\E\oplus\C)$. Denote by $\overline X$ the closure of $X$ in $\P(\E\oplus\C)$. This is a projective variety and $X$ is a dense Zariski open subset of $\overline X$. 
Assume that the claim is not true. Then, there is a sequence of points $z_{[1]}, z_{[2]},\ldots$ in $X$ converging to some point $a\in \overline X\backslash X$ such that  $\pi(z_{[j]})$ converges to some point in $\T$. 

\begin{lemma} \label{l:approx-surface}
There is  a holomorphic map $\tau:\D\to \overline X$ such that $\tau(0)=a$, $\tau(\D^*)\subset X$ and  $\pi(\tau(x))$ admits cluster values in $\T$ when $x$ goes to $0$. Moreover, $\tau(\D^*)$ is Zariski dense in $X$.
\end{lemma}
\proof
Recall that $\Gamma_\R$ is the real space spanned by $\Gamma$. For $t>0$, denote by $\Gamma_\R(t)$ the set of points in $\E$ of distance less than $t$ from $\Gamma_\R$. Observe that we can fix a $t$ large enough so that $\Gamma_\R(t)$ contains $z_{[j]}$ for every $j$. 
Near the point $a$, using a suitable coordinate system, we can see $\Gamma_\R(t)$ as an open cone at $a$ whose intersection with $\overline X$ contains the sequence $z_{[j]}$ which converges to $a$. By the curve selecting lemma  in \cite{Lojasiewicz},  there is a real analytic curve $\tau:[-1,1]\to \overline X$ such that 
$\tau(0)=a$, $\tau([-1,1]\backslash \{0\}) \subset X$ and $\tau((0,1])\subset \Gamma(2t)$. 

Since $\tau$ is real analytic, we can extend it to a holomorphic map from a neighbourhood of $0\in\C$ to $\overline X$. Without lost of generality, we can assume that $\tau$ is defined on $\D$. Observe that the set $\tau^{-1}(\overline X\backslash X)$ should be discrete in $\D$. Therefore, by reducing the size of $\D$, we can assume that $\tau(\D^*)\subset X$. It is not difficult to see that $\pi(\tau((0,1]))$ belongs to $\pi(\Gamma_\R(2t))$ which is relatively compact in $\T$. So $\tau$ satisfies the first part of the lemma.

Now, it is enough to replace $\tau$ by a   map from $\D$ to $\overline X$ which is equal to $\tau$ at 0 up to a high enough order. We still have $\tau(\D^*)\subset X$, $\tau(0)=a$, $\tau((0,1]) \subset \Gamma_\R(3t)$ and moreover
 the second property in the lemma holds for a choice of $\tau$, using suitable transcendental holomorphic functions. This proves the lemma.
\endproof

By Lemma \ref{l:approx-surface}, we have that  $\overline{\pi(X)}^\Zariski$ is also the Zariski closure of $\pi(\tau(\D^*))$ in $\T$. Since 
 $\pi(\tau(x))$ admits cluster values in $\T$ when $x$ goes to $0$, the smallest complex subspace $\F$ of $\E$ such that $\tau_\F: \D^* \to \E/ \F$ can be extended holomorphically through $0$ is non-trivial.  By Theorem \ref{t:ALW-1D}, $\overline{\pi(X)}^\Zariski$ is invariant by $\overline{\pi(\F)}^\Zariski.$  On the other hand, by Lemma \ref{l:inv-Lie} and the comment right after Theorem \ref{t:1D-Gamma},  $\overline{\pi(\F)}^\Zariski$ is in fact a  non-trivial sub-semi-torus of $\T.$  This contradicts our assumption at the beginning of the proof of Theorem \ref{t:ALW}. The result follows.
\endproof

\end{document}